\documentclass{elsart}

\usepackage{amssymb,amsmath}

\let\noi=\noindent
\def\init{\mbox{in}}

\begin{document}

\begin{frontmatter}



\title{Valuations in algebraic field extensions}
\author[FJH]{F. J. Herrera Govantes\thanksref{MTM}\thanksref{HF}}
\ead{jherrera@algebra.us.es},
\author[FJH]{M. A. Olalla Acosta\thanksref{MTM}\thanksref{HF}}
\ead{miguelolalla@algebra.us.es},
\author[MS]{M. Spivakovsky\thanksref{HF}}
\ead{spivakov@picard.ups-tlse.fr}
\thanks[MTM]{Partially supported by MTM2004-07203-C02-01 and FEDER}
\thanks[HF]{Partially supported by HF2004-0117.}

\address[FJH]{Departamento de {\'A}lgebra, Facultad de
Matem{\'a}ticas, Universidad de Sevilla, Apartado 1160, 41080 Sevilla,
Spain}


\address[MS]{Laboratoire de Math{\'e}matiques E. Picard, UMR 5580 du
CNRS, Universit{\'e} Paul Sabatier, 31062 Toulouse cedex 9, France}




\begin{abstract}
Let $K\rightarrow L$ be an algebraic field extension and $\nu$ a valuation of $K$. The purpose of this paper is to
describe the totality of extensions $\left\{\nu'\right\}$ of $\nu$ to $L$ using a refined version of MacLane's key
polynomials. In the basic case when $L$ is a finite separable extension and $rk\ \nu=1$, we give an explicit description of
the limit key polynomials (which can be viewed as a generalization of the Artin--Schreier polynomials). We also give a
realistic upper bound on the order type of the set of key polynomials. Namely, we show that if $char\ K=0$ then the set of
key polynomials has order type at most $\mathbb N$, while in the case $char\ K=p>0$ this order type is bounded above by
$\left(\left[\log_pn\right]+1\right)\omega$, where $n=[L:K]$. Our results provide a new point of view of the the well known
formula $\sum\limits_{j=1}^se_jf_jd_j=n$ and the notion of defect.
\end{abstract}

\end{frontmatter}




\section{Introduction}
\label{In}

All the rings in this paper will be commutative with 1.

This paper grew out of the authors' joint work \cite{HOST} with
B. Teissier, which is devoted to classifying the extensions
$\hat\nu$ of a given valuation $\nu$, centered in a local domain
$R$, to rings of the form $\frac{\hat R}P$, where $\hat R$
is the formal completion of $R$ and $P$ is a prime ideal of
$\hat R$ such that $P\cap R=(0)$. In particular, in \cite{HOST}
we are interested in characterizing situations in which the
valuation $\hat\nu$ (or one of its composed valuations) is
unique. This naturally led us to the following

\noindent{\bf Question.} Given a field extension $K\hookrightarrow L$ and a finite rank
valuation $\nu$ of $K$, when is the extension of $\nu$ to $L$ unique?

An obvious necessary condition for uniqueness is that $L$ be algebraic over $K$.

In the present paper, we give an algorithm for describing the totality of extensions $\nu'$
of $\nu$ to $L$ in terms of (a refined version of) MacLane's key polynomials, assuming $L$ is
algebraic over $K$. The case of purely inseparable extensions
being trivial, we will assume that $L$ is separable over $K$.
Since an arbitrary algebraic field extension is a direct limit of
finite extensions, we may assume that $L$ is finite over $K$. In
particular, $L$ is simple by the primitive element theorem; write $L=K[x]$.

It is sufficient to solve the problem in the case $rk\ \nu=1$: the case of a valuation of an
arbitrary finite rank will then follow by induction on $rk\ \nu$. Indeed, if $\nu$ is the
composition of two lower rank valuations $\nu_1$ and $\nu_2$, then $\nu'$ is the composition
of $\nu'_1$ and $\nu'_2$, where $\nu'_1$ is an extension of $\nu_1$ to $L$ and $\nu'_2$ is a
valuation of the residue field $k_{\nu'_1}$ of the valuation ring $R_{\nu'_1}$,
extending $\nu_2$. Since the field $k_{\nu'_1}$ is an algebraic extension of
$k_{\nu_1}$ (\cite{ZS}, Chapter VI, \S11), it is enough to describe the extensions
$\nu'_1$ of $\nu_1$ to $L$ and the extensions of $\nu_2$ to $k_{\nu'_1}$.

Two main techniques used in this paper are higher Newton polygons and a version of MacLane's
key polynomials, similar to those considered by M. Vaqui\'e (\cite{V1}, \cite{V2}, \cite{V3}, and \cite{V4}), and
reminiscent of related objects studied by Abhyankar and Moh (approximate roots \cite{AM1}, \cite{AM2}) and T.C. Kuo
(\cite{Kuo1}, \cite{Kuo2}). When $L=K[x]$
is a simple extension of $K$, our algorithm is phrased in terms of the slopes of higher Newton polygons of
the minimal polynomial $f$ of $x$, the first one being the usual Newton polygon of $x$; the
algorithm amounts to successively constructing key polynomials of $\nu'$. At each step of the
algorithm there are finitely many possibilities to choose from. Namely, at the $i$-th step we
have to choose a non-vertical side $L$ of the $i$-th Newton polygon, consider the polynomial
$g$ over the graded algebra of $\nu$ determined by $L$ and choose an irreducible factor of
$g$. The number of steps itself can be countable (in fact, the number of steps has
order type at most $\omega$ in characteristic zero and is bounded above by the ordinal
$(\left[\log_pn\right]+1)\omega$ in characteristic $p>0$, where $n$ is the degree of $x$ over
$K$ and $\omega$ stands for the first infinite ordinal). Thus our algorithm can be viewed as providing an
answer to the above question about uniqueness: the extension $\nu'$ is unique
if and only if the choice of both $L$ and $g$ is unique at every step of the
algorithm. A simple sufficient condition for the extension $\nu'$ to be unique
is that the image $\init_{\nu'}x$ in the graded algebra $G_{\nu'}$ have the same degree $n$ over the graded
algebra $G_\nu$ of $\nu$ as $x$ does over $K$; this condition is valid whether
or not $\nu$ has rank 1 and has a very explicit characterizaiton in terms of
the (first) Newton polygon of $f$ (namely, it is equivalent to saying that the
Newton polygon has only one non-vertical side $L$ and the polynomial over the
graded algebra of $\nu$, determined by $L$, is irreducible).

This paper is organized as follows. In \S\ref{algebras} we summarize some
basic definitions and results about algebras without zero divisors, graded by
ordered semigroups. \S\ref{key}--\S\ref{positive} are devoted to the main
construction of the paper --- that of key polynomials. Namely, we suppose
given an extension $\nu'$ of $\nu$ to $L$. We define a well ordered set
$\mathbf Q=\{Q_i\}_{i\in\Lambda}$ of key polynomials of $\nu'$, which may be
finite or countable. If $char\ K=0$, the set $\Lambda$ has order type at most
$\omega$; if $char\ K=p>0$ then $\Lambda$ has order type strictly less than
$(\left[\log_pn\right]+1)\omega$, where $n$ is the degree of $x$ over $K$.

\noi\textbf{Notation.} $\mathbb N$ will denote the set of non-negative integers.
For an element $l\in\Lambda$, we will denote by $l+1$
the immediate successor of $l$ in $\Lambda$. The immediate predecessor, when
it exists, will be denoted by $l-1$. For a positive integer $t$, $l+t$ will
denote the immediate successor of $l+(t-1)$. For an element $l\in\Lambda$, the
initial segment $\{Q_i\}_{i<l}$ of the set of key polynomials will be denoted
by $\mathbf Q_l$. Throughout this paper, we let
\begin{eqnarray}
p&=&1\qquad\qquad\ \ \,\text{if }char\ K=0\\
&=&char\ K\qquad\text{if }char\ K>0.
\end{eqnarray}
In \S\ref{key}, we will fix an ordinal $l$ and assume that the key polynomials $\mathbf Q_{l+1}$ are already defined. We
will define the notion of the $l$-th Newton polygon and the $l$-standard expansion of an element of $K[X]$ with respect to
$\mathbf Q_{l+1}$. We will then define the next key polynomial $Q_{l+1}$. Roughly speaking, $Q_{l+1}$ will be defined to be
the lifting to $L$ of the monic minimal polynomial, satisfied by $\init_{\nu'}Q_l$
over the graded algebra $G_\nu\left[\init_{\nu'}\mathbf Q_l\right]$.

In \S\ref{infinite} we study the situation when the
above recursive algorithm does not stop after finitely many steps, that is, when it gives rise to an infinite sequence
$\{Q_{l+t}\}_{t\in\mathbb N}$ of key polynomials. We define a pair $(\delta_i(f),\epsilon_i(f))$ of basic positive integer
invariants of the Newton polygon $\Delta_i(f)$ (where $i$ runs over
the set of all ordinals for which $Q_i$ is defined). We prove that the pair $(\delta_i(f),\epsilon_i(f))$ is non-increasing
in the lexicographical ordering. We deduce that if $char\ K=0$ and $rk\ \nu=1$ then iterating this construction
at most $\omega$ times, we obtain a sequence $\{Q_i\}$ of key polynomials such that
\begin{equation}
\lim\limits_{i\to\infty}\nu'(Q_i)=\infty.\label{eq:infty}
\end{equation}
In \S\ref{diff} we study the effect of the differential operators $\frac1{p^b!}\frac{\partial^{p^b}}{\partial x^{p^b}}$
on key polynomials and on $f$ in the case the above invariant $\delta_i(f)$ stabilizes.

In \S\ref{tendstoinfty} we use the results of \S5 to show that $\delta_i(f)$ can stabilize only if it is of the form
$\delta_i(f)=p^e$ for some $e\in\mathbb N$.

In \S\ref{positive} we assume that $char\ K=p>0$ and consider an ordinal $l$
which does not have an immediate predecessor. We assume that the key
polynomials $\mathbf Q_l$ are already defined and then define the next key
polynomial $Q_l$. We show that this case can occur at most $[\log_pn]$
times. A set of key polynomials is said to be \textbf{complete} if every $\nu'$-ideal of
$R_{\nu'}$ is generated by products of powers of the $Q_i$ (in other words,
the valuation $\nu'$ is completely determined by the data $\{Q_i,\nu'(Q_i)\}$).
In \S\ref{keyproperty} we prove the main property of key polynomials
$\{Q_i\}$, constructed in \S\S\ref{key}--\ref{positive}: they form a complete set of key
polynomials.

An algorithm for describing the totality of extensions $\nu'$ can be read off
from this data. This algorithm will be described in \S\ref{descr}. As a corollary, we deduce the well known formula
$\sum\limits_{j=1}^se_if_id_i=n$, where $\{\nu_1,\dots,\nu_s\}$ is the set of all the extensions of $\nu$ to $L$, $f_i$
is the index of the value group of $\nu$ viewed as a subgroup of the value group of $\nu_i$, $e_i$ is the degree of the
reside field extension $k_\nu\hookrightarrow k_{\nu_i}$ and $d_i$ is the defect of $\nu_i$ (a much more complete and
detailed treatment of this formula can be found in M. Vaqui\'e's paper \cite{V4}).

In case $char\ K=p>0$ our algorithm is less satisfactory than in characteristic zero
in that at certain junctures it depends on non-constructive considerations
such as a given subset of $\Gamma$ having a maximum or an upper bound.

The idea of using key polynomials and Newton polygons in this context is not
new. What we believe to be new in this paper is the explicit description of
the totality of key polynomials and the definition and an explicit
construction of limit key polynomials, rather intricate in the case of
positive characteristic. In particular, we believe that our bound on the order
type of the set of key polynomials required is new and is the first realistic
bound of its kind.

We want to acknowledge the fact that there is some intersection of our results with
those obtained independently and simultaneously by Michel Vaqui\'e \cite{V4}. We thank him for helpful conversations and, in
particular, for sharing his insights into the notion of defect.

\section{Algebras graded by ordered semigroups.}
\label{algebras}

Graded algebras associated to valuations will play a crucial role in this
paper. In this section, we give some basic definitions and prove several easy results about graded
algebras. Throughout this paper, a ``graded algebra'' will mean ``an
algebra without zero divisors, graded by an ordered semigroup''. As
usual, for a graded algebra $G$, $ord$ will denote the natural valuation of
$G$, given by the grading.
\begin{defn} Let $G$ be a graded algebra without zero
divisors. The {\bf saturation} of $G$, denoted by $G^*$, is the
graded algebra
$$
G^*=\left\{\left.\frac gh\ \right| \ g,h\in G,\
h\text{ homogeneous},\ h\ne0\right\}.
$$
The algebra $G$ is said to be {\bf saturated} if $G=G^*$.
\end{defn}
Of course, we have $G^*=(G^*)^*$ for any graded algebra $G$, so
$G^*$ is always saturated.

The main example of saturated graded algebras appearing in this paper is the
following.
\begin{exmp}Let $\nu:K^*\rightarrow\Gamma$ be a valuation. Let
$(R_\nu,M_\nu,k_\nu)$ denote the valuation ring of $\nu$. For $\beta\in\Gamma$, consider the
following $R_\nu$-submodules of $K$:
$$
\aligned
\mathbf P_\beta&=\{x\in K^*\ |\ \nu(x)\ge\beta\}\cup\{0\},\\
\mathbf P_{\beta+}&=\{x\in K^*\ |\ \nu(x)>\beta\}\cup\{0\}.
\endaligned
$$
We define
$$
G_\nu=\bigoplus_{\beta\in\Gamma}\frac{\mathbf P_\beta}{\mathbf P_{\beta+}}.
$$
The $k_\nu$-algebra $G_\nu$ is an integral domain. For any element $x\in K^*$
with $\nu(x)=\beta$, the natural image of $x$ in $\frac{\mathbf P_\beta}{\mathbf P_{\beta+}}\subset G_\nu$ is a
homogeneous element of $G_\nu$ of degree $\beta$, which we will denote by
$\init_\nu x$. The algebra $G_\nu$ is saturated.

Let $\nu'$ be an extension of $\nu$ to $L$. For an element $\beta\in\Gamma$,
let
\begin{eqnarray}
\mathbf P'_\beta&=&\left\{y\in L\ \left|\
\nu'(x)\ge\beta\right.\right\}\cup\{0\}\\
\mathbf P'_{\beta+}&=&\left\{y\in L\ \left|\
\nu'(x)>\beta\right.\right\}\cup\{0\}.
\end{eqnarray}
Put $G_{\nu'}=\bigoplus\limits_{\beta\in\Gamma}\frac{\mathbf P'_\beta}{\mathbf
P'_{\beta+}}$.
The extension $G_\nu\rightarrow G_{\nu'}$ of graded algebras is finite of
degree bounded by $[L:K]$ (cf. \cite{ZS}, Chapter VI, \S11). In the present paper, we do not
use this result of Zariski--Samuel but rather give another proof of it.
\end{exmp}
\begin{rem} Let $G,G'$ be two graded algebras without zero
divisors, with $G\subset G'$. Let $x$ be a homogeneous element of
$G'$, satisfying an algebraic dependence relation
\begin{equation}
a_0x^\alpha+a_1x^{\alpha-1}+\dots+a_\alpha=0\label{eq:relation}
\end{equation}
over $G$ (here $a_j\in G$ for $0\le j\le\alpha$). Without loss of
generality, we may assume that (\ref{eq:relation}) is homogeneous (that is, the
quantity $j\ ord\ x+ord\ a_j$ is constant for $0\le j\le\alpha$; this is
achieved by replacing (\ref{eq:relation}) by the sum of those terms $a_jx^j$
for which the quantity $j\ ord\ x+ord\ a_j$ is minimal), and that the integer
$\alpha$ is the smallest possible. Dividing (\ref{eq:relation}) by $a_0$, we
see that $x$ satisfies an {\it integral} homogeneous relation over $G^*$ of
degree $\alpha$ and no algebraic relation of degree less than $\alpha$. In
other words, $x$ is {\it algebraic} over $G$ if and only if it is {\it
integral} over $G^*$; the conditions of being ``algebraic over $G^*$'' and
``integral over $G^*$'' are one and the same thing.
\end{rem}
Let $G\subset G'$, let $x$ be as above and let $G[x]$ denote the graded
subalgebra of $G'$, generated by $x$ over $G$. By the above Remark, we may
assume that $x$ satisfies a homogeneous integral relation
\begin{equation}
x^\alpha+a_1x^{\alpha-1}+\dots+a_\alpha=0\label{eq:relation1}
\end{equation}
over $G^*$ and no algebraic relations over $G^*$ of degree less than $\alpha$.
\begin{prop}\label{inverse} Every element of $(G[x])^*$ can be written uniquely as a polynomial
in $x$ with coefficients in $G^*$, of degree strictly less than $\alpha$.
\end{prop}
\begin{pf} Let $y$ be a homogeneous element of $G[x]$. Since $x$ is integral over
$G^*$, so is $y$. Let
\begin{equation}
y^\gamma+b_1y^{\gamma-1}+\dots+b_\gamma=0\label{eq:relation2}
\end{equation}
with $b_j\in G^*$, be a homogeneous integral dependence relation of $y$ over $G^*$, with
$b_\gamma\ne0$. By (\ref{eq:relation2}),
$$
\frac1y=-\frac1{b_\gamma}(y^{\gamma-1}+b_1y^{\gamma-2}+\dots+b_{\gamma-1}).
$$
Thus, for any $z\in G[x]$, we have
\begin{equation}
\frac zy\in G^*[x].\label{eq:frac}
\end{equation}
Since $y$ was an arbitrary homogeneous element of $G[x]$, we have
proved that
$$
(G[x])^*=G^*[x].
$$
Now, for every element $y\in G^*[x]$ we can add a multiple of (\ref{eq:relation1})
to $y$ so as to express $y$ as a polynomial in $x$ of degree less than
$\alpha$. Moreover, this expression is unique because $x$ does not
satisfy any algebraic relation over $G^*$ of degree less than $\alpha$.
\qed
\end{pf}

The following result is an immediate consequence of definitions:
\begin{prop}\label{nondeg} Let $G_\nu$ be the graded algebra associated to a valuation
$\nu:K\rightarrow\Gamma$, as above. Consider a sum of the form $y=\sum\limits_{i=1}^sy_i$,
with $y_i\in K$. Let $\beta=\min\limits_{1\le i\le s}\nu(y_i)$ and
$$
S=\left\{\left.i\in\{1,\dots,n\}\ \right|\ \nu(y_i)=\beta\right\}.
$$
The following two conditions are equivalent:

(1) $\nu(y)=\beta$

(2) $\sum\limits_{i\in S}\init_\nu y_i\ne0$.
\end{prop}

\section{Key polynomials and higher Newton polygons}
\label{key}

Let $K\hookrightarrow L$ be a finite separable field extension and
$\nu:K^*\rightarrow\Gamma$ a valuation of $K$ of real rank 1, where $\Gamma$ is a $\mathbb
Q$-divisible group and $\nu(K^*)$ is a subgroup of $\Gamma$. The extension $L$ is simple
by  the primitive element theorem. Pick and fix a generator $x$ of $L$ over $K$; write
$L=K[x]$.

In this section we begin the main construction of the paper --- that of key polynomials.
Namely, we suppose given an extension $\nu'$ of $\nu$ to $L$.

\begin{defn} A \textbf{complete set of key polynomials} for $\nu'$ is a well ordered
collection $\mathbf Q=\{Q_i\}_{i\in\Lambda}$ of elements of $L$ such that for each
$\beta\in\Gamma$ the $R_\nu$-module $\mathbf P'_\beta$ is generated by all the products of
the form $\prod\limits_{j=1}^sQ_{i_j}^{\gamma_j}$ such that
$\sum\limits_{j=1}^s\gamma_j\nu'(Q_{i_j})\ge\beta$.
\end{defn}
Note, in particular, that if $\mathbf Q$ is a complete set of key polynomials
then their images $\init_{\nu'}Q_i\in G_{\nu'}$ induce a set of generators of $G_{\nu'}$
over $G_\nu$. Furthermore, we want to make the set $\Lambda$ as small as possible, that is, to minimize the order type of
$\Lambda$.

Our algorithm for constructing all the possible extensions $\nu'$ of $\nu$ to $L$ amounts
to successively constructing key polynomials until the resulting set of key polynomials
becomes complete for $\nu'$.

We will fix an ordinal $l$ and assume that the key polynomials
$\mathbf Q_{l+1}$ are already defined (the notation $\mathbf Q_{l+1}$ is defined in the
Introduction). We will then define the next key polynomial $Q_{l+1}$. If $Q_{l+1}=0$, the algorithm stops. In
\S\ref{infinite} we will study what happens when this algorithm does not stop after finitely many steps and
will show that if $char\ K=0$ then iterating this construction at most $\omega$ times,
we obtain a sequence $\{Q_i\}_{i\in\mathbb N}$ of elements of $L$ such that
\begin{equation}
\lim\limits_{i\to\infty}\nu'(Q_i)=\infty.\label{eq:infty1}
\end{equation}
This will end the construction of key polynomials in characteristic zero; in
\S\ref{keyproperty} we will show that the resulting set of key polynomials is complete.

For each $l\in\Lambda$, we will define the notion of the $l$-th Newton polygon
and the $l$-standard expansion of an element of $K[x]$ with respect to $\mathbf Q_{l+1}$.
Roughly speaking, $Q_{l+1}$ will be defined to be the lifting to $L$ of the monic minimal
polynomial, satisfied by $\init_{\nu'}Q_l$ over the graded algebra
$G_\nu\left[\init_{\nu'}\mathbf Q_l\right]$. An algorithm for describing the totality of
extensions $\nu'$ can be read off from this data. This algorithm will be described in
\S\ref{descr}.

Put $Q_1=x$ and $\alpha_1=1$.

Let $X$ be an independent variable and let $f=\sum\limits_{i=0}^na_iX^i$
denote the minimal polynomial of $x$ over $K$. Making a change of variables of
the form $x\rightarrow ax$ with $a\in K$, if necessary, we may assume that
\begin{equation}
\nu(a_n)<\nu(a_i)\text{ for }0\le i<n;\label{eq:normaliz}
\end{equation}
furthermore, dividing $f$ by $a_n$ we may assume $f$ to be monic with
$\nu(a_i)>0$ for $0\le i<n$. The condition (\ref{eq:normaliz}) is needed to
ensure that
\begin{equation}
\nu'(x)>0\label{eq:positive}
\end{equation}
for any extension $\nu'$ of $\nu$ to $L$. Let $\Gamma_+$ (resp. $\mathbb Q_+$)
denote the semigroup of non-negative elements of $\Gamma$ (resp. $\mathbb Q$).

Take an element $h=\sum\limits_{i=0}^sd_iX^i\in K[X]$.
\begin{defn} The first \textbf{Newton polygon} of $h$ with respect to $\nu$ is
the convex hull $\Delta_1(h)$ of the set
$\bigcup\limits_{i=0}^s\left(\left(\nu(d_i),i\right)+
\left(\Gamma_+\oplus\mathbb Q_+\right)\right)$ in $\Gamma\oplus\mathbb Q$.
\end{defn}
To an element $\beta_1\in\Gamma_+$, we associate the following valuation
$\nu_1$ of $K(X)$: for a polynomial $h=\sum\limits_{i=0}^sd_iX^i$, we put
$$
\nu_1(h)=\min\left\{\left.\nu(d_i)+i\beta_1\ \right|\ 0\le i\le s\right\}.
$$
In what follows, for an element $y\in L$, we will write informally $\nu_1(y)$
for $\nu_1(y(X))$, where $y(X)$ is the unique representative of $y$ in
$K[X]$ of degree strictly less than $n$. Similarly, for a polynomial $h\in
K[X]$ we will sometimes write $\nu'(h)$ to mean $\nu'(h\mod\ (f))$.

Consider an element $\beta_1\in\Gamma_+$.
\begin{defn} We say that $\beta_1$ \textbf{determines a side} of $\Delta_1(h)$
if the following condition holds. Let
$$
S_1(h,\beta_1)=\left\{\left.i\in\{0,\dots,s\}\ \right|\
i\beta_1+\nu(d_i)=\nu_1(h)\right\}.
$$
We require that $\#S_1(h,\beta_1)\ge2$.
\end{defn}
Let $\beta_1=\nu'(x)$. Then for any $h\in K[X]$ we have $\nu_1(h)\le\nu'(h)$;
furthermore, $\nu_1(f)<\infty=\nu'(f)$.
\begin{prop}\label{vanish} Take a polynomial $h=\sum\limits_{i=0}^sd_iX^i\in
K[X]$ such that
\begin{equation}
\nu_1(h)<\nu'(h)\label{eq:strict}
\end{equation}
(for example, we may take $h=f$). Then
\begin{equation}
\sum\limits_{i\in S(h,\beta_1)}\init_\nu d_i\init_{\nu'}x^i=0.
\end{equation}
\end{prop}
\begin{pf} We have
$$\sum\limits_{i\in S(h,\beta_1)}d_ix^i=h(x)-\sum\limits_{i\in\{0,\dots,s\}\setminus
S(h,\beta_1)}d_ix^i,$$
hence
$$\nu'\left(\sum\limits_{i\in
S(h,\beta_1)}d_ix^i\right)>\nu_1(h).$$
Then $\sum\limits_{i\in
S_1(h,\beta_1)}\init_\nu d_i\init_{\nu'}x^i=0$ in $\frac{\mathbf
P'_{\nu_1(h)}}{\mathbf P'_{\nu_1(h)}+}\subset G_{\nu'}$ by
Proposition
\ref{nondeg}.
\qed
\end{pf}
\begin{cor} Take a polynomial $h\in K[X]$ such that $\nu_1(h)<\nu '(h)$. Then
$\beta_1$ determines a side of $\Delta_1(h)$.
\end{cor}
\begin{pf} If $S_1(h,\beta_1)$ consisted of a single element $i$, we would have
$$\init_\nu d_i\init_{\nu'}x^i\ne0,$$
contradicting Proposition \ref{vanish}.
\qed
\end{pf}

Letting $h=f$, we see from (\ref{eq:normaliz}) that
$\beta_1>0$ (geometrically, this corresponds to the fact that the side of
$\Delta_1(f)$ determined by $\beta_1$ has strictly negative slope).

\noindent\textbf{Notation.} Let $\bar X$ be a new variable. Take a polynomial
$h$ as above. We denote
$$
\init_1h:=\sum\limits_{i\in S_1(h,\beta_1)}\init_\nu d_i\bar X^i.
$$
The polynomial $\init_1h$ is quasi-homogeneous in $G_\nu[\bar X]$, where the
weight assigned to $\bar X$ is $\beta_1$. Let
\begin{equation}
\init_1h=v\prod\limits_{j=1}^tg_j^{\gamma_j}\label{eq:factoriz}
\end{equation}
be the factorization of $\init_1h$ into irreducible factors in $G_\nu[\bar
X]$. Here $v\in G_\nu$ and the $g_j$ are monic polynomials in $G_\nu[\bar X]$
(to be precise, we first factor $\init_1h$ over the field of fractions of
$G_\nu$ and then observe that all the factors are quasi-homogeneous and
therefore lie in $G_\nu[\bar X]$).
\begin{prop}\label{minimal} (1) The element $\init_{\nu'}x$ is integral over
$G_\nu$.

(2) The minimal polynomial of $\init_{\nu'}x$ over $G_\nu$ is one of the
    irreducible factors $g_j$ of (\ref{eq:factoriz}).

\end{prop}
\begin{pf} Both (1) and (2) of the Proposition follow from the fact
that $\init_{\nu'}x$ is a root of the polynomial $\init_1h$ (Proposition
\ref{vanish}).
\qed
\end{pf}
Now take $h=f$. Renumbering the factors in (\ref{eq:factoriz}), if necessary,
we may assume that $g_1$ is the minimal polynomial of $\init_{\nu'}x$ over
$G_\nu$. Let $\alpha_2=\deg_{\bar X}g_1$. Write
$g_1=\sum\limits_{i=0}^{\alpha_2}\bar b_i\bar X^i$, where $\bar
b_{\alpha_2}=1$. For each $i$, $0\le i\le\alpha_2$, choose a representative
$b_i$ of $\bar b_i$ in $R_\nu$ (that is, an element of $R_\nu$ such that
$\init_\nu b_i=\bar b_i$; in particular, we take $b_{\alpha_2}=1$). Put
$Q_2=\sum\limits_{i=0}^{\alpha_2}b_ix^i$.
\begin{defn} The elements $Q_1$ and $Q_2$ are called, respectively,
\textbf{the first and second key polynomials} of $\nu'$.
\end{defn}
Now, every element $y$ of $L$ can be written uniquely as a finite sum of the
form
\begin{equation}
y=\sum_{\begin{array}{c}0\le\gamma_1<\alpha_2\\
\gamma_1+\gamma_2\alpha_2<n\end{array}}
b_{\gamma_1\gamma_2}Q_1^{\gamma_1}Q_2^{\gamma_2}\label{eq:standard2}
\end{equation}
where $b_{\gamma_1\gamma_2}\in K$ (this is proved by Euclidean division by the
monic polynomial $Q_2$). The expression (\ref{eq:standard2}) is called
\textbf{the second standard expansion of }$y$.

Now, take an ordinal number $l\ge2$ which has an immediate predecessor; denote
this ordinal by $l+1$. If $char\ K=0$, assume that $l\in\mathbb N$. Assume,
inductively, that key polynomials $\mathbf Q_{l+1}$, and positive integers
$\mbox{\boldmath$\alpha$}_{l+1}=\{\alpha_i\}_{i\le l}$ are already
constructed, and that all but finitely many of the $\alpha_i$ are equal to
1. We want to define the key polynomial $Q_{l+1}$.

We will use the following multi-index notation:
$\mbox{\boldmath$\gamma$}_{l+1}=\{\gamma_i\}_{i\le l}$, where all but finitely
many $\gamma_i$ are equal to 0, $\mathbf
Q_{l+1}^{\mbox{\boldmath$\gamma$}_{l+1}}=\prod\limits_{i\le
l}Q_i^{\gamma_i}$. Let $\beta_i=\nu'(Q_i)$.
\begin{defn} An index $i<l$ is said to be $l$-\textbf{essential} if there
exists a positive integer $t$ such that either $i+t=l$ or $i+t<l$ and
$\alpha_{i+t}>1$; otherwise $i$ is called $l$-\textbf{inessential}.
\end{defn}
In other words, $i$ is $l$-inessential if and only if $i+\omega\le l$ and
$\alpha_{i+t}=1$ for all $t\in\mathbb N$.

\noi\textbf{Notation.} For $i<l$, let
\begin{eqnarray}
i+&=&i+1\qquad\,\text{if }i\text{ is $l$-essential}\\
&=&i+\omega\qquad\text{otherwise}.
\end{eqnarray}
\begin{defn}\label{lstandard} A multiindex $\mbox{\boldmath$\gamma$}_{l+1}$ is
said to be \textbf{standard with respect to}
$\mbox{\boldmath$\alpha$}_{l+1}$ if
\begin{equation}
0\le\gamma_i<\alpha_{i+}\text{ for }i\le l,\label{eq:standind0}
\end{equation}
\begin{equation}
\sum\limits_{i\le l}\gamma_i\prod\limits_{j\le i}\alpha_j\le n,\label{eq:standind}
\end{equation}
and if $i$ is $l$-inessential then the set $\{j<i+\ |\ j+=i+\text{ and }\gamma_j\ne0\}$ has cardinality at most one.
An $l$\textbf{-standard monomial in} $\mathbf Q_{l+1}$ (resp. an
$l$\textbf{-standard monomial in} $\init_{\nu'}\mathbf Q_{l+1}$) is a product
of the form $c_{\mbox{\boldmath$\gamma$}_{l+1}}\mathbf
Q_{l+1}^{\mbox{\boldmath$\gamma$}_{l+1}}$,
(resp. $c_{\mbox{\boldmath$\gamma$}_{l+1}}\init_{\nu'}\mathbf
Q_{l+1}^{\mbox{\boldmath$\gamma$}_{l+1}}$) where
$c_{\mbox{\boldmath$\gamma$}_{l+1}}\in K$
(resp. $c_{\mbox{\boldmath$\gamma$}_{l+1}}\in G_\nu$) and the multiindex
$\mbox{\boldmath$\gamma$}_{l+1}$ is standard with respect to
$\mbox{\boldmath$\alpha$}_{l+1}$.
\end{defn}
\begin{rem} In the case when $i$ admits an immediate predecessor, the condition
(\ref{eq:standind0}) amounts to saying that $\gamma_{i-1}<\alpha_i$.
\end{rem}
\begin{defn}\label{notinvolving1} An $l$\textbf{-standard expansion not involving }$Q_l$ is a
finite sum $S$ of $l$-standard monomials, not involving $Q_l$, having the
following property. Write $S=\sum\limits_\beta S_\beta$, where $\beta$ ranges
over a certain finite subset of $\Gamma_+$ and
\begin{equation}
S_\beta=\sum\limits_jd_{\beta j}\label{eq:notinvolving}
\end{equation}
is a sum of standard monomials $d_{\beta j}$ of value $\beta$. We require that
\begin{equation}
\sum\limits_j\init_{\nu'}d_{\beta j}\ne0\label{eq:nonzero}
\end{equation}
for each $\beta$ appearing in (\ref{eq:notinvolving}).
\end{defn}
In the special case when $l\in\mathbb N$, (\ref{eq:nonzero}) holds
automatically for any sum of $l$-standard monomials not involving $Q_l$ (this
follows from Proposition \ref{value} below by induction on $l$).
\begin{prop}\label{notinvolving} Let $l$ be an ordinal and $t$ a positive
integer. Assume that the key polynomials $\mathbf Q_{l+t+1}$ are defined and
that $\alpha_l=\dots=\alpha_{l+t}=1$. Then any $(l+t)$-standard expansion
does not involve any $Q_i$ with $l\le i<l+t$. In particular, an $l$-standard expansion not involving $Q_l$ is the same
thing as an $(l+t)$-standard expansion, not involving $Q_{l+t}$.
\end{prop}
\begin{pf} (\ref{eq:standind0}) implies that for $l\le i\le l+t$,
$Q_i$ cannot appear in an $(l+t)$-standard expansion with a positive
exponent. \qed
\end{pf}
We will frequently use this fact in the sequel without mentioning it explicitly.
\begin{defn}\label{standard} For an element $g\in K[X]$, an expression of the
form $g=\sum\limits_{j=0}^sc_jQ_l^j$, where each $c_j$ is an $l$-standard
expansion not involving $Q_l$, will be called an $l$-\textbf{standard
expansion of }$g$. For a non-zero element $y\in L$, an $l$-\textbf{standard
expansion of }$y$ is an $l$-standard expansion of the representative $y(X)$ of
$y$ in $K[X]$ of degree strictly less than $n$.
\end{defn}
In what follows, we will be mostly interested in standard expansions of
non-zero elements of $L$ and of the polynomial $f(X)$.
\begin{defn} Let $\sum\limits_{\mbox{\boldmath$\gamma$}}\bar
  c_{\mbox{\boldmath$\gamma$}}\init_{\nu'}\mathbf
  Q_{l+1}^{\mbox{\boldmath$\gamma$}}$ be an $l$-standard expansion, where
  $\bar c_{\mbox{\boldmath$\gamma$}}\in G_\nu$. A \textbf{lifting} of
  $\sum\limits_{\mbox{\boldmath$\gamma$}}\bar
  c_{\mbox{\boldmath$\gamma$}}\init_{\nu'}\mathbf
  Q_{l+1}^{\mbox{\boldmath$\gamma$}}$ to $L$ is an $l$-standard expansion
  $\sum\limits_{\mbox{\boldmath$\gamma$}}c_{\mbox{\boldmath$\gamma$}}\mathbf
  Q_{l+1}^{\mbox{\boldmath$\gamma$}}$, where $c_{\mbox{\boldmath$\gamma$}}$ is
  a representative of $\bar c_{\mbox{\boldmath$\gamma$}}$ in $K$.
\end{defn}
\begin{defn} Assume that $char\ K=p>0$. An $l$-standard expansion\linebreak
  $\sum\limits_jc_jQ_l^j$, where each $c_j$ is an $l$-standard expansion not
  involving $Q_l$, is said to be \textbf{weakly affine} if $c_j=0$ whenever
  $j>0$ and $j$ is not of the form $p^e$ for some $e\in\mathbb N$.
\end{defn}
Assume, inductively, that for each ordinal $i\le l$, every element $h$ of $L$ and
the polynomial $f(X)$ admit an $i$-standard expansion. Furthermore, assume that for each $i\le l$, the $i$-th key polynomial
$Q_i$ admits an $i$-standard expansion, having the following additional properties.

If $i$ has an immediate predecessor $i-1$ in $\Lambda$ (such is always the
case in characteristic 0), the $i$-th standard expansion of $Q_i$ has the form
\begin{equation}
Q_i=Q_{i-1}^{\alpha_i}+\sum\limits_{j=0}^{\alpha_i-1}
\left(\sum\limits_{\mbox{\boldmath$\gamma$}_{i-1}}
c_{ji\mbox{\boldmath$\gamma$}_{i-1}}\mathbf
Q_{i-1}^{\mbox{\boldmath$\gamma$}_{i-1}}\right)Q_{i-1}^j,
\label{eq:standform}
\end{equation}
where:

\noi(1) each $c_{ji\mbox{\boldmath$\gamma$}_{i-1}}\mathbf
Q_{i-1}^{\mbox{\boldmath$\gamma$}_{i-1}}$ is an $(i-1)$-standard monomial, not involving $Q_{i-1}$

\noi(2) the quantity $j\beta_{i-1}+\sum\limits_{q<i-1}\gamma_q\beta_q$ is constant
for all the monomials
$$\left(c_{ji\mbox{\boldmath$\gamma$}_{i-1}}\mathbf
Q_{i-1}^{\mbox{\boldmath$\gamma$}_{i-1}}\right)Q_{i-1}^j$$
appearing on the right hand side of (\ref{eq:standform})

\noi(3) the equation
\begin{equation}
\init_{\nu'}Q_{i-1}^{\alpha_i}+\sum\limits_{j=0}^{\alpha_i-1}
\left(\sum\limits_{\mbox{\boldmath$\gamma$}_{i-1}}\init_\nu
c_{ji\mbox{\boldmath$\gamma$}_{i-1}}\init_{\nu'}\mathbf
Q_{i-1}^{\mbox{\boldmath$\gamma$}_{i-1}}\right)
\init_{\nu'}Q_{i-1}^j=0\label{eq:mini-1}
\end{equation}
is the minimal algebraic relation satisfied by $\init_{\nu'}Q_{i-1}$ over the
subalgebra $G_\nu[\init_{\nu'}\mathbf Q_{i-1}]^*\subset G_{\nu'}$.

Finally, if $char\ K=p>0$ and $i$ does not have an immediate predecessor in
$\Lambda$ then there exist an $i$-inessential index $i_0$ and a strictly
positive integer $e_i$ such that $i=i_0+$ and
$Q_i=\sum\limits_{j=0}^{e_i}c_{ji_0}Q_{i_0}^{p^j}$ is a weakly affine monic
$i_0$-standard expansion of degree $\alpha_i=p^{e_i}$ in $Q_{i_0}$, where each
$c_{ji_0}$ is an $i_0$-standard expansion not involving $Q_{i_0}$. Moreover,
there exists a positive element $\bar\beta_i\in\Gamma$ such that
\begin{eqnarray}
\bar\beta_i&>&\beta_q\qquad\text{ for all }q<i,\label{eq:Qi1}\\
\beta_i&\ge&p^{e_i}\bar\beta_i\quad\text{ and}\label{eq:Qi2}\\
p^j\bar\beta_i+\nu(c_{i_0j})&=&p^{e_i}\bar\beta_i\quad\text{ for }0\le j\le
e_i.\label{eq:Qi3}
\end{eqnarray}
If $i\in\mathbb N$, we assume, inductively, that the $i$-standard expansion is unique. If $char\ K>0$, and
$h=\sum\limits_{j=0}^{s_i}d_{ji}Q_i^j$ is an $i$-standard expansion of $h$ (where $h$ is either $f(X)$ or an element of
$L$), we assume that the elements $d_{ji}\in L$ are uniquely determined by $h$ (strictly speaking, this
does not mean that the $i$-standard expansion is unique: for example, if $i$ is a limit ordinal, $d_{ji}$ admits an
$i_0$-standard expansion for each $i_0<i$ such that $i=i_0+$, but there may be countably many choices of
$i_0$ for which such an $i_0$-standard expansion is an $i_0$-standard expansion, not
involving $Q_{i-1}$ in the sense of Definition \ref{notinvolving1}).
\begin{prop}\label{bound} (1) The polynomial $Q_i$ is monic in $x$; we have
\begin{equation}
\deg_xQ_i=\prod\limits_{j\le i}\alpha_j.\label{eq:degQi}
\end{equation}
(2) Let $z$ be an $i$-standard expansion, not involving $Q_i$. Then
\begin{equation}
\deg_xz<\deg_xQ_i.\label{eq:degreeless}
\end{equation}
\end{prop}
\begin{pf} (\ref{eq:degQi}) and (\ref{eq:degreeless}) are proved simultaneously by
transfinite induction on $i$, using (\ref{eq:standform}) and (\ref{eq:standind0}) repeatedly
to calculate and bound the degree in $x$ of any standard monomial (recall that by assumption
all but finitely many of the $\alpha_i$ are equal to 1).
\qed
\end{pf}
The rest of this section is devoted to the definition of $Q_{l+1}$. In what
follows, we will sometimes not distinguish between the elements $Q_i$ and
their representatives in $K[X]$ in order to simplify the notation. When we do wish to make
such a distinction, we will denote the representative of $Q_i$ in $K[X]$ by $Q_i(X)$.

Write
\begin{equation}
f=\sum\limits_{i=0}^{n_l}a_{jl}Q_l^j,\label{eq:flstandard}
\end{equation}
where each $a_{jl}$ is a homogeneous $l$-standard expansion not involving $Q_l$, such that
$$
\deg_xa_{jl}+j\prod\limits_{q=1}^l\alpha_q\le n,
$$
with strict inequality for $j<n_l$.

Take any element $h\in K[X]$ and let $h=\sum\limits_{i=0}^sd_iQ_l^i$ be an
$l$-standard expansion of $h$, where each $d_i$ is an $l$-standard expansion,
not involving $Q_l$.
\begin{defn} The $l$-th Newton polygon of $h$ with respect to $\nu$ is the
convex hull $\Delta_l(h)$ of the set
$\bigcup\limits_{i=0}^s\left(\left(\nu'(d_i),i\right)+
\left(\Gamma_+\oplus\mathbb Q_+\right)\right)$ in $\Gamma\oplus\mathbb Q$.
\end{defn}
To an element $\beta_l\in\Gamma_+$, we associate a valuation $\nu_l$ of $K(X)$
as follows. Given an $l$-standard expansion $h=\sum\limits_{i=0}^sd_iQ_l^i$ as
above, put $\nu_l(h)=\min\limits_{0\le i\le s}\{i\beta_l+\nu'(d_i)\}$. Note
that even though in the case of positive characteristic the standard
expansions of the elements $d_i$ are not, in general, unique, the elements
$d_i\in L$ themselves are unique by Euclidean division, so $\nu_l$ is well
defined. That $\nu_l$ is, in fact, a valuation, rather than a
pseudo-valuation, follows from the definition of standard expansion,
particularly, from (\ref{eq:nonzero}). We always have $\nu_l(h)\le\nu'(h)$ and
$\nu_l(f)<\infty=\nu'(f)$.

\noindent\textbf{Notation.} Let $\bar Q_l$ be a new variable and let $h$ be as
above. We denote
\begin{eqnarray}
S_l(h,\beta_l):&=&\left\{j\in\{0,\dots,s\}\ \left|\
j\beta_l+\nu'(d_j)=\nu_l(h)\right.\right\}.\label{eq:Ll}\\
\init_lh:&=&\sum\limits_{j\in S_l(h,\beta_l)}\init_{\nu'}d_j\bar
Q_l^j;\label{eq:initl}
\end{eqnarray}
The polynomial $\init_lh$ is quasi-homogeneous in
$G\left[\init_{\nu'}\mathbf Q_l,\bar Q_l\right]$, where the weight assigned
to $\bar Q_l$ is $\beta_l$.

Take a polynomial $h$ such that
\begin{equation}
\nu_l(h)<\nu'(h)\label{eq:greaterl}
\end{equation}
(for example, we may take $h=f$).
\begin{prop}\label{vanishl} We have $\sum\limits_{j\in
S_l(h,\beta_l)}\init_{\nu'}\left(d_jQ_l^j\right)=0$ in $\frac{\mathbf
P'_{\nu_l(h)}}{\mathbf P'_{\nu_l(h)+}}\subset G_{\nu'}$.
\end{prop}
\begin{pf} This follows immediately from (\ref{eq:greaterl}), the
fact that
$$
\sum\limits_{j\in S_l(h,\beta_l)}d_jQ_l^j=h-\sum\limits_{j\in
S_l(h,\beta_l)\setminus\{0,\dots,s\}}d_jQ_l^j
$$
and Proposition \ref{nondeg}.
\qed
\end{pf}
Let $\beta_l$ be a non-negative element of $\Gamma$.
\begin{defn} We say that $\beta_l$ \textbf{determines a side} of $\Delta_l(h)$
if $\#S_l(h,\beta_l)\ge2$.
\end{defn}
\begin{cor}\label{betaldetermines} Let $\beta_l=\nu'(Q_l)$. Then:

(1) $\beta_l$ determines a side of $\Delta_l(h)$.

(2)
\begin{eqnarray}
\beta_l&>&\alpha_l\beta_{l-1}\quad\text{ if }\ (l-1)\text{ exists}\\
\beta_l&\ge&p^{e_l}\bar\beta_l\qquad\text{otherwise}.
\end{eqnarray}
\end{cor}
\begin{pf} (1) Suppose not. Then the sum
$0=\sum\limits_{j\in S_l(h,\beta_l)}\init_{\nu'}\left(d_jQ_l^j\right)$
consists of only one term and hence cannot be 0. This contradicts Proposition
\ref{vanishl}; (1) is proved.

(2) follows immediately from (\ref{eq:mini-1}) and (\ref{eq:Qi2}). This
completes the proof of Corollary \ref{betaldetermines}.
\qed
\end{pf}
Let
\begin{equation}
\init_lh=v_l\prod\limits_{j=1}^tg_{jl}^{\gamma_{jl}}\label{eq:factorizl}
\end{equation}
be the factorization of $\init_lh$ into (monic) irreducible
factors in $G_\nu\left[\init_{\nu'}\mathbf Q_l\right]\left[\bar
Q_l\right]$ (to be precise, we first factor $\init_lh$ over the
field of fractions of $G_\nu\left[\init_{\nu'}\mathbf Q_l\right]$
and then observe that all the factors are quasi-homogeneous and
therefore lie in $G_\nu\left[\init_{\nu'}\mathbf Q_l\right]\left[\bar Q_l\right]$).
\begin{cor} The element $\init_{\nu'}Q_l$ is integral over $G_\nu$. Its
minimal polynomial over $G_\nu$ is one of the irreducible factors $g_{jl}$ of
(\ref{eq:factorizl}).
\end{cor}
Put $h=f$ in (\ref{eq:factorizl}). Renumbering the factors in
(\ref{eq:factorizl}), if necessary, we may assume that $g_{1l}$ is the minimal
polynomial of $\init_{\nu'}Q_l$ over $G_\nu\left[\init_{\nu'}\mathbf
Q_l\right]$. Let
\begin{equation}
\alpha_{l+1}=\deg_{\bar Q_l}g_{1l}.\label{eq:alphal+1}
\end{equation}
Write
\begin{equation}
g_{1l}=\bar Q_l^{\alpha_{l+1}}+\sum\limits_{j=0}^{\alpha_{l+1}-1}
\left(\sum\limits_{\mbox{\boldmath$\gamma$}_l}\bar
c_{l+1,j\mbox{\boldmath$\gamma$}_l}\init_{\nu'}\mathbf
Q_l^{\mbox{\boldmath$\gamma$}_l}\right) \bar Q_l^j,\label{eq:ingraded}
\end{equation}
Define \textbf{the $(l+1)$-st key polynomial} of $\nu'$ to be a lifting
\begin{equation}
Q_{l+1}=Q_l^{\alpha_{l+1}}+\sum\limits_{j=0}^{\alpha_{l+1}-1}
\left(\sum\limits_{\mbox{\boldmath$\gamma$}_l}
c_{l+1,j\mbox{\boldmath$\gamma$}_l}\mathbf
Q_l^{\mbox{\boldmath$\gamma$}_l}\right)Q_l^j\label{eq:lifting}
\end{equation}
(\ref{eq:ingraded}) to $L$. In the special case when $t=\alpha_{l+1}=1$ in (\ref{eq:factorizl}) and (\ref{eq:alphal+1}),
some additional (and rather intricate) conditions must be imposed on the lifting (\ref{eq:lifting}). In fact, in this case
we will define several consecutive key polynomials at the same time. We will now explain what these additional conditions
are, after making one general remark:
\begin{rem} Since $g_{1l}$ is an irreducible polynomial in $\bar Q_l$ by definition, the key polynomial $Q_{l+1}(X)$ is
also irreducible (for a non-trivial factorization of $Q_{l+1}(X)$ would give rise to a non-trivial factorization of
$g_{1l}$).
\end{rem}
To define $Q_{l+1}$ in the case $t=\alpha_{l+1}=1$, we first introduce two numerical
characters of the situation which will play a crucial role in the rest of the paper. Let
$\delta_l(h)=\deg_{\bar Q_l}\init_lh$.
\begin{defn} The vertex $\left(\nu'\left(a_{\delta_l(h),l}\right),\delta_l(h)\right)$ of the
Newton polygon $\Delta(h)$ is called the \textbf{pivotal vertex} of $\Delta(h)$.
\end{defn}
Let
\begin{equation}
\nu_l^+(h)=\min\left\{\left.\nu'\left(d_jQ_l^j\right)\ \right|\ \delta_l(h)<j\le
s\right\}\label{eq:nu+}
\end{equation}
and
$$
S'_l(h)=\left\{j\in\{\delta_l(h)+1,\dots,s\}\ \left|\
\nu'\left(d_jQ_l^j\right)=\nu_l^+(h)\right.\right\}.
$$
Let $\epsilon_l(h)=\max\ S'_l(h)$ (if the set on the right hand side of (\ref{eq:nu+}) is
empty, we adopt the convention that $\nu_l^+(h)=\epsilon_l(h)=\infty$). The quantities
$\delta_l(h)$ and $\epsilon_l(h)$ are strictly positive by definition. It follows from
definitions that $\epsilon_l(h)>\delta_l(h)$. Below, we will see that that the pair
$(\delta_l(h),\epsilon_l(h))$ is non-increasing with $l$ (in the lexicographical ordering),
that the equality $\delta_{l+1}(h)=\delta_l(h)$ imposes strong restrictions on $\init_lh$ and
that decreasing $(\delta_l(f),\epsilon_l(h))$ strictly ensures that the algorithm stops
after a finite number of steps.

Assume that $t=\alpha_{l+1}=1$ in (\ref{eq:factorizl}) and (\ref{eq:alphal+1}). Let
$\delta=\delta_l(f)$. We have $v_l=\init_{\nu'}a_{\delta l}$ and (\ref{eq:factorizl})
rewrites as
\begin{equation}
\init_lf=\init_{\nu'}a_{\delta l}g_{1l}^\delta.\label{eq:factoriz2}
\end{equation}
In what follows, we will consider $l$-standard expansions of the form
\begin{equation}
Q'=Q_l+z_l+\dots+z_i,\label{eq:tildez}
\end{equation}
where each $z_j$ is a homogeneous $l$-standard expansion, not involving $Q_l$, such that
\begin{equation}
\beta_l=\nu'(z_l)<\nu'(z_{l+1})<\dots<\nu'(z_i).\label{eq:inequalities}
\end{equation}
\begin{rem}\label{degree} Note that by (\ref{eq:degreeless}), we have $\deg_xz_q<\deg_xQ_l$
for all $q$.
\end{rem}
\begin{defn} Let $Q'$ be as above. A \textbf{standard expansion} of $f$ with respect to $Q'$
is an expression of the form $f=\sum\limits_{j=0}^{n_l}a'_j{Q'}^j$, where each $a'_j$ is an
$l$-standard expansion, not involving $Q_l$. The \textbf{Newton polygon} $\Delta(f,Q')$ of
$f$ with respect to $Q'$ is the convex hull in $\Gamma_+\oplus\mathbb Q_+$ of the set
$\bigcup\limits_{i=0}^{n_l}\left(\left(\nu'(a'_i),i\right)+
\left(\Gamma_+\oplus\mathbb Q_+\right)\right)$.
\end{defn}
Substituting $Q'-z_l-\dots-z_i$ for $Q_l$ in (\ref{eq:flstandard}), writing
$$
\sum\limits_{i=0}^{n_l}a_{il}(Q'-z_l-\dots-z_i)^i=\sum\limits_{j=0}^{n_l}a'_j{Q'}^j,
$$
and using (\ref{eq:inequalities}), we see that $\nu'(a'_\delta)=\nu'(a_{\delta l})$ and that
$\left(\nu'\left(a_{\delta l}\right),\delta\right)$ is a vertex of
$\Delta(f,Q')$ (though it might not be the pivotal one).
\begin{defn} The \textbf{characteristic side} of $\Delta(f,Q')$ is the side $A(f,Q')$ whose upper endpoint is
$(\nu'(a_{\delta l}),\delta)$.
\end{defn}
Let $\beta(Q')$ denote the element of $\Gamma_+$ which determines the side $A(f,Q')$.

For $l\le j\le i$, put $Q'_j=Q_l+z_l+\dots+z_{j-1}$, let $\Delta(f,Q'_j)$ be the corresponding
Newton polygon and $A(f,Q'_j)$ the characteristic side of $\Delta(f,Q'_j)$.

Let $T$ denote the set of all the $l$-standard expansions of the form (\ref{eq:tildez}),
where each $z_j$ is a homogeneous $l$-standard expansion, not involving $Q_l$, such that
the inequalities (\ref{eq:inequalities}) hold, $\nu'(z_i)<\beta(Q')$ and
\begin{equation}
\init_{A\left(f,Q'_j\right)}f=\init_{\nu'}a'_\delta\left(\bar
Q+\init_{\nu'}z_j\right)^\delta\label{eq:inA}
\end{equation}
whenever $l\le j<i$.

We impose the following partial ordering on $T$. Given an element $Q'=Q_l+z_l+\dots+z_i\in
T$ with $i>l$, we declare its immediate predecessor in $T$ to be the element
$Q_l+z_l+\dots+z_{i-1}$. By definition, our partial ordering is the coarsest
one among those in which $Q_l+z_l+\dots+z_{i-1}$ precedes $Q_l+z_l+\dots+z_i$ for all the elements $Q'$ as above.

Take an element $Q':=Q_l+z_l+\dots+z_i\in T$. Let $A'=A(f,Q')$.
\begin{rem}\label{slope} Assume that
$$
\init_{A'}f=\init_{\nu'}a'_\delta(\bar Q+\init_{\nu'}z')^\delta
$$
for some $l$-standard expansion $z'$, not involving $Q_l$. Then
\begin{equation}
\init_{\nu'}(Q')=-\init_{\nu'}z';\label{eq:init}
\end{equation}
in particular, $\nu'(Q')=\nu'(z')$. In other words, $\nu'(Q')$, $\nu'(z')$ and the slope of
the side $A'$ are all equivalent sets of data. In the sequel, we prefer to talk about
$\nu'(z')$ rather than $\nu'(Q')$ for the following reason. In \S\ref{descr}, rather than
working with a fixed valuation $\nu'$, we will use the same algorithm to construct all the
possible extensions $\nu'$. Therefore it will be important to describe the next step in the
algorithm using only the data known at this stage of the construction, rather than the entire
data of $\nu'$ itself. Since we are assuming that the key polynomials $\mathbf Q_l$ and their
values are already known, we may consider $\nu'(z')$ as being known as well.
\end{rem}
\noi\textbf{Notation.} In what follows, for an element $b\in L$, $b(X)$ will denote the
representative of $b$ in $K[X]$ of degree less than $n$.
\begin{prop}\label{iota1} Consider an $l$-standard expansion $w$ of the form
$w_l+w_{l+1}+\dots+w_i$, where $w_l$, \dots, $w_i$ are homogeneous $l$-standard expansions
and $w_l$ is an $l$-standard expansions, not involving $Q_l$, such that
$\beta_l=\nu'(w_l)<\dots<\nu'(w_i)$. Fix an element $\beta\in\Gamma_+$,
\begin{equation}
\beta>\beta_l.\label{eq:betabetal}
\end{equation}
Then $w(X)$ can be written in the form
\begin{equation}
w(X)=w_l(X)+\tilde w_{l+1}(X)+\dots+\tilde
w_j(X)+w^\dag(X)((Q_l(X)+w(X)))+\psi_\beta(X),\label{eq:decomp}
\end{equation}
where $\psi_\beta$ is an $l$-standard expansion, $w_l$, $\tilde w_{l+1}$, \dots,
$\tilde w_j$ are homogeneous $l$-standard expansions, not involving $Q_l$, such that
\begin{eqnarray}
\beta_l=\nu'(w_l)<\nu'(\tilde w_{l+1})<\dots&<&\nu'(\tilde w_j),\label{eq:increasing1}\\
\nu_l(\psi_\beta(X))&\ge&\beta,\qquad\text{ and }\\
\nu_l(w^\dag(X))&>&0\label{eq:nugreater2}
\end{eqnarray}
\end{prop}
\begin{pf} Let $\mu=\nu_l(w_{l+1})$. By definitions, the Proposition is true for
$\beta=\mu$. Assume that the Proposition holds for a certain $\beta$. We will show that it
holds for $\beta$ replaced by $\beta+\mu$, and that will complete the proof. Consider an
expression (\ref{eq:decomp}) satisfying (\ref{eq:increasing1})--(\ref{eq:nugreater2}). Write
$\psi_\beta$ in the form $\psi_\beta=\psi_{\beta+\mu}+\tilde\psi$, where
$\nu_l(\psi_{\beta+\mu}(X))\ge\beta+\mu$ and $\tilde\psi$ consists of monomials of value
greater than or equal to $\beta$ but strictly less than $\beta+\mu$. By assumptions and
Remark \ref{degree}, $\deg_Xw_l(X)<\deg_XQ_l(X)$. Divide the polynomial $\tilde\psi(X)$ by $Q_l(X)+w_l(X)$:
$$
\tilde\psi(X)=q(X)(Q_l(X)+w_l(X))+r(X),
$$
where $\deg_Xr(X)<\deg_XQ_l(X)$. Then
$$
\tilde\psi(X)=q(X)(Q_l(X)+w(X))+r(X)+\tilde\psi_{\beta+\mu}(X),
$$
where $\nu_l\left(\tilde\psi_{\beta+\mu}\right)\ge\beta+\mu$. Absorb the quotient $q(X)$ into
$w^\dag(X)$ and $\tilde\psi_{\beta+\mu}$ into $\psi_{\beta+\mu}$. Let
\begin{equation}
r(X)=\tilde w_{j+1}(X)+\dots+\tilde w_{\tilde j}(X)\label{eq:remainder}
\end{equation}
be the $l$-standard expansion of $r(X)$. Since the remainder $r(X)$ is of degree strictly
less than $\deg_XQ_l(X)$, its standard expansion (\ref{eq:remainder}) does not involve any
monomials divisible by $Q_l(X)$. We obtain the desired decomposition
$$\begin{array}{ccl}
w(X) & = & w_l(X)+\tilde w_{l+1}(X)+\dots+\tilde w_{\tilde j}(X)+\\
& & +(w^\dag(X)+q(X))(Q_l(X)+w(X))+\psi_{\beta+\mu}(X).
\end{array}
$$
Condition (\ref{eq:betabetal}) implies that $\nu_l(w^\dag(X)+q(X))>0$, as desired.
\qed
\end{pf}

\begin{prop}\label{compare} Consider two elements
$Q':=Q_l+z'_l+\dots+z'_{i'},Q'':=Q_l+z''_l+\dots+z''_{i''}\in T$. Let
$\Delta'(f)$ and $\Delta''(f)$ be the corresponding Newton polyhedra and $A'$
(resp. $A''$) the characteristic side of $\Delta'(f)$ (resp. $\Delta''(f)$). Assume that
\begin{eqnarray}
\init_{A'}f&=&\init_{\nu'}a_{\delta l}(\bar Q+\init_{\nu'}w')^\delta\qquad\text{ and}\\
\init_{A''}f&=&\init_{\nu'}a_{\delta l}(\bar Q+\init_{\nu'}w'')^\delta
\end{eqnarray}
for some $l$-standard expansions $w'$ and $w''$, not involving
$Q_l$. Furthermore, assume that
$$
\nu'(Q')<\nu'(Q'').
$$
Then there exists a third element
$$
Q''':=Q_l+z'''_l+\dots+z'''_{i'''}\in T,\ Q'''>Q',
$$
having the following property. Let $\Delta'''(f)$ denote the Newton polygon
determined by $Q'''$ and $A'''$ the characteristic side of $\Delta'''(f)$. Then $A'''=A''$ and
$\init_{A'''}f=\init_{A''}f$.
\end{prop}

\begin{pf} Let $w=Q''-Q'$ and fix an element $\beta\in\Gamma$, $\beta>\nu'(Q'')$.
Apply Proposition \ref{iota1} with $Q_l$  replaced by $Q'$. The hypotheses of Proposition
\ref{iota1} are satisfied because
$$
\nu'(w')=\nu'(Q')<\nu'(Q'')=\nu'(w'')
$$
and $\init_{\nu'}w'=-\init_{\nu'}Q'$ by assumptions, hence
$\init_{\nu'}w=-\init_{\nu'}Q'=\init_{\nu'}w'$, in particular, $\init_{\nu'}w$ does not
involve $\init_{\nu'}Q_l$. By Proposition \ref{iota1} we can write
\begin{equation}
w=z_{i'+1}+\dots+z_{i'''}+w^\dag Q''+\psi_\beta\label{eq:w}
\end{equation}
such that
\begin{eqnarray}
\nu_l(w^\dag)&>&0,\label{eq:wdag}\\
\nu_l(\psi_\beta)&>&\nu'(Q'')\label{eq:psibeta}
\end{eqnarray}
and $z_{i'+1}$, \dots, $z_{i'''}$ are $l$-standard expansions, not involving $Q_l$. Put
$$
Q'''=Q_l+z'_l+\dots+z'_{i'}+z_{i'+1}+\dots+z_{i'''}.
$$
Then (\ref{eq:w}), (\ref{eq:wdag}) and (\ref{eq:psibeta}) show that
$$
\nu'(Q''')=\nu'(Q'')
$$
and
$$
\init_{\nu'}Q'''=\init_{\nu'}Q'';
$$
the Proposition follows immediately.
\qed
\end{pf}

To define $Q_{l+1}$ in the special case when
\begin{equation}
t=\alpha_{l+1}=1\label{eq:talphal+1}
\end{equation}
in (\ref{eq:factorizl}) and (\ref{eq:alphal+1}), first assume that $char\ K=0$. Equations
(\ref{eq:factoriz2}) and (\ref{eq:talphal+1}) imply that $a_{\delta-1,l}\ne0$ and
\begin{equation}
g_{1l}=\bar Q_l+\init_{\nu'}\frac{a_{\delta-1,l}}{\delta\ a_{\delta l}}.\label{eq:linear}
\end{equation}
Consider the $l$-standard expansion of $\frac{a_{\delta-1,l}}{\delta
a_{\delta l}}$ and write it in the form
\begin{equation}
\frac{a_{\delta-1,l}}{\delta\
a_{\delta l}}=z_l+z_{l+1}+\dots+z_{l_1-1}+\phi+w,\label{eq:sequence}
\end{equation}
where $l_1$ is an integer strictly greater than $l$, each $z_i$ is a homogeneous
$l$-standard expansion, not involving $Q_l$, such that
\begin{equation}
\nu'\left(\frac{a_{\delta-1,l}}{\delta\
a_{\delta l}}\right)=\nu'(z_l)<\nu'(z_{l+1})<\dots<\nu'(z_{l_1-1})<\nu_l^+(h)-\nu_l(h)+
\beta_l,\label{eq:<bi}
\end{equation}
$\phi$ is a sum of standard monomials of value greater than or equal to
$\nu_l^+(h)-\nu_l(h)+\beta_l$ and $w$ is divisible by $Q_l+\frac{a_{\delta-1,l}}{\delta
a_{\delta l}}$ (such an expression (\ref{eq:sequence}) exists by Proposition \ref{iota1}). Let
\begin{equation}
s=\max\{i\ |\ l\le i\le l_1\text{ and }Q_l+z_l+\dots+z_{i-1}\in T\}.\label{eq:maxj}
\end{equation}
For $l\le i\le s$, define $Q_i=Q_l+z_l\dots+z_{i-1}$.

Next, assume $char\ K=p>0$. Two cases are possible:

\noi\textbf{Case 1.} The set $T$ contains a maximal element. Let $z=z_l+z_{l+1}+\dots+z_{s-1}$
be this maximal element, where each $z_i$ is a homogeneous $l$-standard expansion, not
involving $Q_l$, and $s$ is an ordinal of the form $s=l+t$, $t\in\mathbb N$. Define
$$
Q_i=Q_l+z_l+\dots+z_{i-1}\qquad\text{ for }l+1\le i\le s.
$$
\textbf{Case 2.} The set $T$ does not contain a maximal element. Let
$$\bar\beta=\sup\{\nu(Q')\ |\ Q'\in T\}$$
(here we allow the possibility $\bar\beta=\infty$).
In this case, Proposition \ref{compare} (together with Remark \ref{slope}) shows that there
exists an infinite sequence $z_l,z_{l+1},\dots$ of homogeneous $l$-standard expansions, not
involving $Q_l$, such that for each $t\in\mathbb N$ we have
\begin{equation}
Q_l+z_l+\dots+z_{l+t}\in T\label{eq:inT}
\end{equation}
and $\lim\limits_{t\to\infty}\nu(Q_l+z_l+\dots+z_{l+t})=\bar\beta$; pick and fix one such
sequence. Define
$$
Q_{l+t}=Q_l+z_l+z_{l+1}+\dots+z_{l+t-1}\qquad\text{ for }t\in\mathbb N.
$$
Note that (\ref{eq:inT}), (\ref{eq:inA}) and Remark \ref{slope} imply that the sequence
$\{\nu(Q_l+z_l+\dots+z_{l+t})\}_{t\in\mathbb N}$ is strictly increasing.

For future reference, it will be convenient to distinguish two subcases of
Case 2:

\noi\textbf{Case 2a.} $\bar\beta=\infty$, that is, the sequence
$\{\beta_{l+t}\}_{t\in\mathbb N}$ is unbounded in $\Gamma$. In this case, the
definition of the key polynomials $Q_i$ is complete. In \S\ref{infinite}, we will use differential operators to show that in
this case $\delta$ is necessarily of the form $p^e$ for some $e\in\mathbb N$.

\noi\textbf{Case 2b.} The set $\{\nu(Q')\ |\ Q'\in T\}$ has a least upper bound
$\bar\beta<\infty$ (but no maximum) in $\Gamma$. In this case, we must continue the
construction and define $Q_{l+\omega},Q_{l+\omega+1}$, etc. This will be accomplished in
\S\ref{positive}.
\begin{rem} Note that the definition of $Q_{l+1}$ depends only on the key polynomials
$\mathbf Q_{l+1}$ defined so far, their values $\mbox{\boldmath$\beta$}_{l+1}$ and the
resulting Newton polygons $\Delta_i(f)$, $i\le l$. This will be important in \S\ref{descr}
where we will use the $Q_i$ to \textit{construct} all the possible extensions $\nu'$.
\end{rem}
\begin{prop}{\label{value}} Let $y$ be an element of $L$, represented by a
polynomial in $K[X]$ of degree strictly less than
$\deg_xQ_{l+1}=\prod\limits_{i=0}^{l+1}\alpha_i$. Then $\nu'(y)=\nu_l(y)$.
\end{prop}
\begin{pf} Let $y=\sum\limits_{j=0}^sc_jQ_l^j$ be an $l$-standard
expansion of $y$, where each $c_j$ is an $l$-standard expansion not involving
$Q_l$. Let
$$
S=\left\{j\in\{0,\dots,s\}\ \left|\ \nu'\left(c_j
Q_l^j\right)=\nu_l(y)\right.\right\}.
$$
Let $\bar c_j:=\init_{\nu'}c_j$. Since the degree of $\init_{\nu'}Q_l$ over
$G_\nu[\init_{\nu'}\mathbf Q_l]^*$ is $\alpha_{l+1}$, we see, using the assumption on
$\deg_xy$, that $\sum\limits_{j=0}^s\bar c_j\init_{\nu'}Q_l^j\ne0$ in
$G_{\nu'}$. The result now follows from Proposition \ref{nondeg}.
\qed
\end{pf}

Now, take any polynomial $h\in K[X]$. The $(l+1)$-st standard expansion
$h=\sum\limits_{j=0}^sc_jQ_{l+1}^j$ is constructed from the $l$-th one by
Euclidean division by the polynomial $Q_{l+1}$. Condition
$\nu_l(c_j)=\nu'(c_j)$ required in the definition of standard expansion
(cf. Definition \ref{standard} and (\ref{eq:nonzero})) follows immediately
from the above Proposition and Proposition \ref{bound} (2).

By induction on $t$, this defines key polynomials $Q_{l+t}$ for $t\in\mathbb N$. If for some
$t\in\mathbb N$ we obtain $Q_{l+t}=0$ in $L$, stop. In \S\ref{keyproperty}, we will show that
$\mathbf Q_{l+t}$ is a complete set of key polynomials for $\nu'$, and, in particular, that
the data $\mathbf Q_{l+t}$ and $\mbox{\boldmath$\beta$}_{l+t}$ completely determines $\nu'$.

If $Q_{l+t}\ne0$ for all $t\in\mathbb N$, we obtain an infinite sequence $\{Q_{l+t}\}$ of
key polynomials. If
\begin{eqnarray}
char\ K&=&0\qquad\text{or}\label{eq:char0}\\
\lim\limits_{t\to\infty}\beta_{l+t}&=&\infty,\label{eq:toinfty1}
\end{eqnarray}
stop (in fact, in the next section we will see that (\ref{eq:char0}) implies
(\ref{eq:toinfty1}) and also that in this case $\delta$ has the form $p^e$, $e\in\mathbb N$). In \S\ref{keyproperty}, we
will show that the $\{Q_{l+t}\}$ is a complete set of key polynomials for $\nu'$. If $char K=p>0$ and
$\lim\limits_{t\to\infty}\beta_{l+t}<\infty$, the construction of the next key polynomial
$Q_{l+\omega}$ will be described in \S\ref{positive}.

In the next three sections, we analyze the case when infinitely many such iterations
give rise to an infinite sequence $\{Q_{l+i}\}$ of key polynomials.

\section{Infinite sequences of key polynomials.}
\label{infinite}

Keep the assumption $rk\ \nu=1$. In this section, we analyze the case
when iterating the recursive construction of the previous section produces an
infinite sequence $\{Q_{l+t}\}_{t\in\mathbb N}$. If $char\ K=0$, we show that if the above
algorithm produces an infinite sequence of key polynomials then
\begin{equation}
\lim\limits_{i\to\infty}\beta_i=\infty.\label{eq:unbounded}
\end{equation}
In \S\ref{keyproperty} we will show that (\ref{eq:unbounded}) implies that the
valuation $\nu'$ is completely determined by the resulting data $\{Q_i\}$ and
$\{\beta_i\}$, that is, that the resulting set $\{Q_i\}$ is, indeed, a complete set of
key polynomials. The case when $char\ K=p>0$ and the values $\beta_i$ are bounded above in $\Gamma$ is studied in detail in
\S\ref{positive}.

Take an ordinal $i$ such that $Q_i$ and $Q_{i+1}$ are defined.
Take a polynomial $h$ such that $\nu_i(h)<\nu'(h)$ (for example, we may take
$h=f$). Consider the $i$-th Newton polygon of $h$. Let $S_i(h,\beta_i)$ be as
in (\ref{eq:Ll}). Recall the definition of $\delta_i(h)$:
\begin{equation}
\delta_i(h):=\max\{S_i(h,\beta_i)\}.\label{eq:delta}
\end{equation}
Let $h=\sum\limits_{j=0}^{s_i}d_{ji}Q_i^j$ denote the $i$-standard expansion of $h$, where
each $d_{ji}$ is an $l$-standard expansion, not involving $Q_l$. Recall the
definition (\ref{eq:nu+}) of $\nu_i^+(h)$. The next Proposition shows that the pair
$(\delta_i(h),\epsilon_i(h))$ is non-increasing with $i$ (in the lexicographical ordering)
and that the equality $\delta_{i+1}(h)=\delta_i(h)$ imposes strong restrictions on $\init_ih$.
\begin{prop}\label{stability} (1) We have
\begin{equation}
\alpha_{i+1}\delta_{i+1}(h)\le\delta_i(h).\label{eq:stab}
\end{equation}

(2) If $\delta_{i+1}(h)=\delta_i(h)$ then
\begin{eqnarray}
&\epsilon_{i+1}(h)\le\epsilon_i(h),\label{eq:epsilondrops}&\\
&\init_ih=\init_{\nu'}d_{\delta_i(h)i}\left(\bar Q_i+\init_{\nu'}z_i\right)^{\delta_i(h)},&\label{eq:initih}
\end{eqnarray}
where $z_i$ is some $i$-standard expansion not involving $Q_i$, and $\init_{i+1}h$ contains a monomial of the form
$\init_{\nu'}d_{\delta_i(h)i}\bar Q_{i+1}^{\delta_i(h)}$; in particular,
\begin{equation}
\init_{\nu'}d_{\delta_i(h)i}=\init_{\nu'}d_{\delta_i(h),i+1}\label{eq:initequal}
\end{equation}

(3) If
\begin{equation}
(\delta_i(h),\epsilon_i(h))=(\delta_{i+1}(h),\epsilon_{i+1}(h))\label{eq:pairequal}
\end{equation}
then
\begin{equation}
\init_{\nu'}d_{\epsilon_i(h)i}=\init_{\nu'}d_{\epsilon_i(h),i+1}.\label{eq:initequal1}
\end{equation}
\end{prop}
\begin{pf} We start with three Lemmas. First, consider the $(i+1)$-standard expansion
of $h$:
\begin{equation}
h=\sum\limits_{j=0}^sd_{j,i+1}Q_{i+1}^j,\label{eq:standardi+1}
\end{equation}
where the $d_{j,i+1}$ are $(i+1)$-standard expansions, not involving $Q_{i+1}$.

\begin{lem}\label{nui} (1) We have
$$
\nu_i(h)=\min\limits_{0\le j\le s}\nu_i\left(d_{j,i+1}Q_{i+1}^j\right)=\min\limits_{0\le j\le
s}\{\nu'(d_{j,i+1})+j\alpha_{i+1}\beta_i\}.
$$
(2) Let
$$
S_{i,i+1}=\left\{j\in\{0,\dots,s\}\ \left|\ \nu_i\left(d_{j,i+1}Q_{i+1}^j\right)=\nu_i(h)\right.\right\}
$$
and $j_0=\max\ S_{i,i+1}$. Then $\delta_i(h)=\alpha_{i+1}j_0+\deg_{Q_i}d_{j_0,i+1}$.
\end{lem}
\begin{pf} (1) Provisionally, let
$$
\mu=\min\limits_{0\le j\le s}\nu_i\left(d_{j,i+1}Q_{i+1}^j\right)=\min\limits_{0\le j\le
s}\{\nu'(d_{j,i+1})+j\alpha_{i+1}\beta_i\},
$$
$$
S'_{i,i+1}=\left\{j\in\{0,\dots,s\}\ \left|\
\nu_i\left(d_{j,i+1}Q_{i+1}^j\right)=\mu\right.\right\},
$$
$j'=\max\ S'_{i,i+1}$ and $\delta'=\alpha_{i+1}j'+\deg_{Q_i}d_{j',i+1}$. We want to show that $\mu=\nu_i(h)$,
$S'_{i,i+1}=S_{i,i+1}$, $j'=j_0$ and $\delta_i(h)=\delta'$.

Let $\bar h=\sum\limits_{j\in S'(i,i+1)}d_{j,i+1}Q_{i+1}^j$. Then
$\nu_i(h-\bar h)>\mu$ by definition, so to prove that $\nu_i(h)=\mu$ it is
sufficient to prove that $\nu_i(\bar h)=\mu$.

Now, $\deg_x\bar h=\deg_xd_{\delta',i+1}Q_i^{\delta'}$ by definition of $\delta'$ and
Proposition \ref{bound} (2). Hence the $i$-standard expansion of $\bar h$ contains the
monomial $d_{\delta',i+1}Q_i^{\alpha_{i+1}\delta'}$ and all the other monomials have degree in $x$
strictly smaller than $\deg_xd_{\delta',i+1}Q_i^{\alpha_{i+1}\delta'}$. Thus
$\nu_i\left(\bar h\right)\le\nu_i\left(d_{\delta',i+1}Q_i^{\alpha_{i+1}\delta'}\right)=\mu$, so
$\nu_i(h)\le\mu$. The opposite inequality is trivial and (1) is proved. (2)
follows immediately from this. \qed
\end{pf}
\begin{lem}\label{tejera} Consider two terms of the form $dQ_{i+1}^j$ and $d'Q_{i+1}^{j'}$
(where $j,j'\in\mathbb N$ and $d$ and $d'$ are $i$-standard expansions not involving
$Q_i$). Assume that
\begin{equation}
\nu_i\left(dQ_{i+1}^j\right)\le\nu_i\left(d'Q_{i+1}^{j'}\right)\label{eq:muless}
\end{equation}
and
\begin{equation}
\nu'\left(dQ_{i+1}^j\right)\ge\nu'\left(d'Q_{i+1}^{j'}\right).\label{eq:nugreater}
\end{equation}
Then $j\ge j'$. If at least one of the inequalities (\ref{eq:muless}),(\ref{eq:nugreater})
is strict then $j>j'$.
\end{lem}
\begin{pf} Subtract (\ref{eq:muless}) from
(\ref{eq:nugreater}) and use the definition of $\nu_i$ and the facts that
$\nu_i(Q_{i+1})=\beta_i$ and $\alpha_{i+1}\beta_i<\beta_{i+1}$. \qed
\end{pf}

In the notation of Lemma \ref{nui}, let $\theta_{i+1}(h)=\min\ S_{i,i+1}$.
\begin{defn}\label{characteristic} The vertex $(\nu(d_{\theta_{i+1}(h),i+1}),\theta_{i+1}(h))$ is called the
\textbf{characteristic vertex} of $\Delta_{i+1}(h)$. By convention, $\theta_1(f)=n$, so the characteristic
vertex of $\Delta_1(f)$ is also defined.
\end{defn}
The notion of characteristic vertex will be needed in \S\ref{descr} when we discuss the totality of the extensions $\nu'$
of $\nu$ and the formula $\sum\limits_jf_je_jd_j=n$. It is important that the characteristic vertex of $\Delta_{i+1}(f)$ is
determined by $\mathbf Q_{i+2}$ and $\mbox{\boldmath$\beta$}_{i+1}$: it does not depend on $\beta_{i+1}$.

Let
\begin{equation}
\init_ih=\init_{\nu'}d_{\delta i}\prod\limits_{j=1}^tg_{ji}^{\gamma_{ji}}\label{eq:factoriz3}
\end{equation}
be the factorization of $\init_ih$ into (monic) irreducible
factors in $G_\nu\left[\init_{\nu'}\mathbf Q_i\right]\left[\bar Q_i\right]$, where $g_{1i}$ is the minimal polynomial of
$\init_{\nu'}Q_i$ over $G_\nu\left[\init_{\nu'}\mathbf Q_i\right]$.
\begin{lem}\label{factors} We have
\begin{equation}
\gamma_{1i}=\theta_{i+1}(h)\label{eq:gammatheta}
\end{equation}
(in particular, $d_{\gamma_{1i},i+1}\ne0$) and
\begin{equation}
\init_{\nu'}d_{\theta_{i+1}(h),i+1}=
\init_{\nu'}d_{\delta i}\prod\limits_{j=2}^tg_{ji}^{\gamma_{ji}}(\init_{\nu'}Q_i).\label{eq:initd1}
\end{equation}
\end{lem}
\begin{pf} Write
$$h=\sum\limits_{q\in S_{i,i+1}}d_{q,i+1}Q_{i+1}^q+\sum\limits_{q\in\{0,\dots,n_{i+1}\}\setminus
S_{i,i+1}}d_{q,i+1}Q_{i+1}^q.$$
By Lemma \ref{nui},
\begin{equation}
\init_ih=\sum\limits_{q\in S_{i,i+1}}\init_id_{q,i+1}\init_iQ_{i+1}^q.\label{eq:inith}
\end{equation}
By definition of $\theta_{i+1}(h)$, $\init_iQ_{i+1}^{\theta_{i+1}(h)}$\ \  is\ \ the\ \ highest\ \ power\ \ of\ \
$\init_iQ_{i+1}$\ \ dividing\linebreak$\sum\limits_{q\in S_{i,i+1}}\init_id_{q,i+1}\init_iQ_{i+1}^q$. Also by definition, we
have
\begin{equation}
\init_iQ_{i+1}=g_{1i}.\label{eq:initQi+1}
\end{equation}
Now (\ref{eq:gammatheta}) follows from (\ref{eq:inith}). Also from (\ref{eq:inith}), we see
that $\init_{\nu'}d_{\theta_{i+1}(h),i+1}$ is obtained by substituting $\init_{\nu'}Q_i$ in $\init_ih$, and (\ref{eq:initd1})
follows. \qed
\end{pf}

Now,\ \ apply\ \ Lemma\ \ \ref{tejera}\ \ to\ \ the\ \ monomials\ \ $d_{\theta_{i+1}(h),i+1}Q_{i+1}^{\theta_{i+1}(h)}$\ \
and\linebreak $d_{\delta_{i+1}(h),i+1}Q_{i+1}^{\delta_{i+1}(h)}$. We have
\begin{equation}
\nu'\left(d_{\delta_{i+1}(h),i+1}Q_{i+1}^{\delta_{i+1}(h)}\right)\le
\nu'\left(d_{\theta_{i+1}(h),i+1}Q_{i+1}^{\theta_{i+1}(h)}\right)\label{eq:nu'less}
\end{equation}
by definition of $\delta_{i+1}$ and
\begin{equation}
\nu_i\left(d_{\theta_{i+1}(h),i+1}Q_{i+1}^{\theta_{i+1}(h)}\right)=\nu_i(h)\le
\nu_i\left(d_{\delta_{i+1}(h),i+1}Q_{i+1}^{\delta_{i+1}(h)}\right)\label{eq:nuiless}
\end{equation}
by Lemma \ref{nui}, so the hypotheses of Lemma \ref{tejera} are satisfied. By Lemma \ref{tejera},
\begin{equation}
\theta_{i+1}(h)\ge\delta_{i+1}(h).\label{eq:charpivot}
\end{equation}
Since
\begin{equation}
\alpha_{i+1}\gamma_{1i}=\alpha_{i+1}\theta_{i+1}(h)\le\deg_{\bar Q_i}\init_ih=\delta_i(h)\label{eq:equality}
\end{equation}
by (\ref{eq:factoriz3}), (1) of the Proposition follows.

(2) Assume that $\delta_{i+1}(h)=\delta_i(h)$. Then the above two monomials coincide and
\begin{equation}
\alpha_{i+1}=1.\label{eq:alpha1}
\end{equation}
Furthermore, we have equality in (\ref{eq:equality}), so $\init_ih=\init_{\nu'}d_{\delta_i(h)i}g_{1i}^{\delta_i(h)}$.
Combined with (\ref{eq:alpha1}), this proves (2) of the Proposition.

Finally, (\ref{eq:epsilondrops}) (assuming
(\ref{eq:stab})) is proved by exactly the same reasoning as (\ref{eq:stab}).
(\ref{eq:initequal1}) (assuming (\ref{eq:pairequal})) is proved by the same reasoning as
(\ref{eq:initequal}). This completes the proof of the Proposition. \qed
\end{pf}

\begin{rem} One way of interpreting Lemma \ref{tejera}, together with the inequalities (\ref{eq:nu'less}),\ \
(\ref{eq:nuiless})\ \ and\ \ (\ref{eq:charpivot})\ \ is\ \ to\ \,say\ \,that\ \,the\ \,characteristic\ \, vertex\linebreak
$(\nu'(d_{\theta_{i+1}(h),i+1}),\theta_{i+1}(h))$ of $\Delta_{i+1}(h)$\ \,always\ \,lies\ \,above\ \,its\ \,
pivotal vertex\linebreak $(\nu'(d_{\delta_{i+1}(h),i+1}),\delta_{i+1}(h))$. This fact will be important in \S\ref{descr}.
\end{rem}
For the rest of this section, assume that $\mathbf Q_{l+1}$ is
defined for a certain ordinal number $l$ and that $\mathbb N$ iterations of
the algorithm of the previous section produce an infinite sequence
$\{Q_{l+t}\}_{t\in\mathbb N}$.

Take an ordinal $i$ of the form $l+t$, $t\in\mathbb N$.
\begin{cor}\textbf{(of Proposition \ref{stability})}\label{alpha1} We have $\alpha_{l+i}=1$ for $i\gg0$.
\end{cor}
This fact can also be easily seen without using Proposition \ref{stability}. Indeed, equations (\ref{eq:degQi}),
(\ref{eq:initl}) and (\ref{eq:factorizl}) show that
$$
\prod\limits_{j\le i}\alpha_j\le n
$$
for all $i$. The Corollary follows immediately.
\qed

Choose the ordinal $l$ above so that $\alpha_{l+t}=1$ for all (strictly)
positive integers $t$. By definition, for $t\in\mathbb N$, we have
\begin{equation}
Q_{l+t+1}=Q_{l+t}+z_{l+t},\label{eq:Qi+1}
\end{equation}
where $z_{l+t}$ is a homogeneous $l$-standard expansion of value
$\beta_{l+t}$, not involving $Q_l$ (cf. Proposition \ref{notinvolving}). By Proposition \ref{bound} (2), we have
\begin{equation}
\deg_xz_{l+t}<\deg_xQ_{l+t}.\label{eq:degree}
\end{equation}
Finally,
\begin{equation}
\init_{\nu'}Q_{l+t}=-\init_{\nu'}z_{l+t}\label{eq:Qz}
\end{equation}
by (\ref{eq:init}).

As before, let $h=\sum\limits_{j=0}^{s_i}d_{ji}Q_i^j$ be an $i$-standard
expansion of $h$ for $i\ge l$, where each $d_{ji}$ is an $l$-standard expansion, not
involving $Q_l$. Note that since $\alpha_{l+t}=1$ for $t\in\mathbb N$, we have
$\deg_xQ_i=\prod\limits_{j=2}^{\alpha_i}\alpha_j=
\prod\limits_{j=2}^{\alpha_l}\alpha_j=\deg_xQ_l$ and so
\begin{equation}
s_i=\left[\frac{\deg_xh}{\deg_xQ_i}\right]=\left[\frac{\deg_xh}{\deg_xQ_l}\right]=s_l.\label{eq:degreeequal}
\end{equation}
By Proposition \ref{stability} (1), $\delta_i(h)$ is constant for all $i\gg l$. Let
$\delta=\delta_i(h)$ for $i\gg l$. Write $\delta=p^eu$, where if $p>1$ then $p\not|u$.
Then, according to Proposition \ref{stability} (2) and using the notation of
(\ref{eq:initl}), we see that for $i\gg l$
\begin{equation}
\delta-p^e\in S_i(h,\beta_i)\label{eq:delta-1}
\end{equation}
(in particular, $d_{\delta-p^e,i}\ne0$) and that
\begin{equation}
\init_iz_i=\left(\frac{\init_id_{\delta-p^e,i}}{u\ \init_id_{\delta i}}\right)^{\frac1{p^e}}\label{eq:zi}
\end{equation}
In what follows, the ordinal $i$ will run over the sequence $\{l+t\}_{t\in\mathbb N}$.

Next, we prove a comparison result which expresses the coefficients $d_{ji}$ in terms of $d_{jl}$ for $\delta-p^e\le
j\le\delta$, modulo terms of sufficiently high value.
\begin{prop}\label{diffop} Assume that
\begin{equation}
\delta_{i+1}(h)=\delta_l(h)=\delta.\label{eq:deltastable}
\end{equation}
Take an integer $v\in\{\delta-p^e,\delta-p^e+1,\dots,\delta\}$.
We have
\begin{equation}
\begin{array}{rcl}
d_{vi}&\equiv&\sum\limits_{j=0}^{\delta-v}(-1)^j\binom{v+j}jd_{v+j,l}
(z_l+\dots+z_{i-1})^j\\
& &\mod\ \mathbf P'_{(\nu_l(h)-v\beta_l)+\min\{\nu_l^+(h)-\nu_l(h),\beta_i-\beta_l\}}.\label{eq:closev}
\end{array}
\end{equation}
In particular, letting $v=\delta-p^e$ and $v=\delta$ in (\ref{eq:closev}) we obtain
\begin{equation}
\begin{array}{rcl}
d_{\delta-p^e,i}&\equiv&\sum\limits_{j=0}^{p^e}(-1)^j\binom{\delta-p^e+j}jd_{\delta-p^e+j,l}
(z_l+\dots+z_{i-1})^j\\
& &\mod\ \mathbf P'_{\nu'(d_{\delta-p^e,l})+\min\{\nu_l^+(h)-\nu_l(h),\beta_i-\beta_l\}}.\label{eq:close-1}
\end{array}
\end{equation}
and
\begin{equation}
d_{\delta i}\equiv d_{\delta l}\ \mod\ \mathbf
P'_{\nu'(d_{\delta l})+\min\{\nu_l^+(h)-\nu_l(h),\beta_i-\beta_l\}},\label{eq:closedelta}
\end{equation}
respectively. If $p^e=1$ (in particular, whenever $char\ K=0$), (\ref{eq:close-1}) reduces to
\begin{equation}
d_{\delta-1,i}\equiv d_{\delta-1,l}-\delta\ d_{\delta l}(z_l+\dots+z_{i-1})\ \mod\ \mathbf
P'_{\nu'(d_{\delta-1,l})+\min\{\nu_l^+(h)-\nu_l(h),\beta_i-\beta_l\}}.\label{eq:close2}
\end{equation}
\end{prop}
\begin{pf} By definitions, we have $Q_i=Q_l+z_l+\dots+z_{i-1}$.
First, we will compare the $l$-standard expansion of $h$ with the $i$-standard one. To
this end, we substitute $Q_l=Q_i-z_l-\dots-z_{i-1}$ into the $l$-standard expansion of
$h$. We obtain
\begin{equation}
h=\sum\limits_{j=0}^{s_l}d_{jl}(Q_i-z_l-\dots-z_{i-1})^j=
\sum\limits_{j=0}^{s_l}d_{ji}Q_i^j.\label{eq:fi-1}
\end{equation}
We want to derive information about $\init_ih$ from (\ref{eq:fi-1}). First note that for each
$q\in\{0,\dots,s_l-1\}$ we have $\deg_x\sum\limits_{j=0}^qd_{jl}(Q_i-z_l-\dots-z_{i-1})^j<(q+1)
\deg_xQ_i$. Hence $d_{q+1,i}$ is completely determined by
$d_{q+1,l},d_{q+2,l},\dots,d_{s_ll}$. Next, for $\delta-v<j\le s_l-v$ and $l\le s\le i-1$, note
that
\begin{equation}
\nu'\left(d_{v+j,l}z_s^j\right)\ge j\beta_l+\nu'(d_{v+j,l})\ge
\nu_l^+(h)-v\beta_l,\label{eq:estimdelta}
\end{equation}
so for $\delta-v<j\le s_l-v$ the terms $d_{v+j,l}Q_l^{v+j}$ in (\ref{eq:fi-1}) contribute nothing to
$$
d_{vi}\mod\ \mathbf P'_{(\nu_l(h)-v\beta_l)+\min\{\nu_l^+(h)-\nu_l(h),\beta_i-\beta_l\}}.
$$
Now, the coefficients $d_{vi}$ in (\ref{eq:fi-1}) are obtained from $\sum\limits_{j=0}^{s_l}d_{jl}(Q_i-z_l-\dots-z_{i-1})^j$
by opening the parentheses and then applying Euclidean division by $Q_i$; such a Euclidean division may change
the coefficients $d_{vi}$ by adding terms of value at least
$\nu'(Q_i)-\nu_l(Q_i)=\beta_i-\alpha_{l+1}\beta_l=\beta_i-\beta_l$. Finally, using (\ref{eq:initih}) (which holds thanks to
the hypothesis (\ref{eq:deltastable})) we observe that for $v$ and $j$ as in (\ref{eq:closev}) we have
$\nu'(d_{v+j,l})\ge\nu'(d_{\delta l})+(\delta-v-j)\beta_l=\nu_l(h)-(v+j)\beta_l$. This completes the proof of
(\ref{eq:closev}).

(\ref{eq:close-1}) and (\ref{eq:closedelta}) follow from (\ref{eq:closev}), after observing that
$$
\nu_l(h)=\nu'(d_{\delta l})+\delta\beta_l=\nu'(d_{\delta-p^e,l})+(\delta-p^e)\beta_l
$$
by (\ref{eq:initih}). (\ref{eq:close2}) obtained from (\ref{eq:close-1}) by substituting $p^e=1$. The Proposition is proved.
\qed
\end{pf}
Now let $f=h$ and let $f=\sum\limits_{j=0}^{n_i}a_{ji}Q_i^j$ be the $i$-standard expansion of $f$. We have $n_i=n_l$ (this
is a special case of (\ref{eq:degreeequal})).
\begin{prop}\label{propinfty} Assume that the sequence $\{Q_i\}$ is infinite. There are
two mutually exclusive possibilities: either
\begin{equation}
\lim\limits_{i\to\infty}\beta_i=\infty\label{eq:infinity}
\end{equation}
or $char\ K=p>0$ and there exists $t_0\in\mathbb N$ such that, letting $i_0=l+t_0$, we have
\begin{equation}
\lim\limits_{\overset\longrightarrow i}\beta_i<\beta_{i_0}+
\frac1{p^e}\left(\nu_{i_0}^+(f)-\nu_{i_0}(f)\right)\label{eq:limbounded}
\end{equation}
(recall that we are assuming $rk\ \nu=1$).
\end{prop}
\begin{pf} We start with a few lemmas.
\begin{lem}\label{deltadrops} Assume that either $char\ K=0$ and $s<l_1$ in (\ref{eq:maxj})
or the set $T$ contains a maximal element $Q_s=Q_l+z_l+\dots+z_{s-1}$. Then
$\delta_{s+1}(f)<\delta$. In particular, this case can occur at most finitely many times.
\end{lem}
\begin{pf} We give a proof by contradiction. Suppose
$\delta_{s+1}(f)=\delta$. By Proposition \ref{stability} (2),
\begin{equation}
\init_sf=\init_{\nu'}a_{\delta s}(\bar Q_s+\init_{\nu'}w)^\delta,\label{eq:wdelta}
\end{equation}
for some $l$-standard expansion $w$, not involving $Q_l$. This shows that $Q_s$ is not maximal
in $T$: the element $Q_s+w$ is greater than $Q_s$. It remains to consider the case
$char\ K=0$ and $s<l_1$. In this case, (\ref{eq:sequence}), (\ref{eq:closedelta}) and
(\ref{eq:close2}) imply that
\begin{eqnarray}
\init_{\nu'}a_{\delta s}&=&\init_{\nu'}a_{\delta l}\qquad\qquad\text{ and}\label{eq:adelta}\\
\init_{\nu'}a_{\delta-1,s}&=&\init_{\nu'}(a_{\delta-1,l}-\delta a_{\delta_l}(z_l+\dots+z_{s-1}))=
\init_{\nu'}(\delta a_\delta z_s)\label{eq:adelta-1}.
\end{eqnarray}
Combining this with (\ref{eq:wdelta}), we see that
$\init_{\nu'}w=\init_{\nu'}z_s$. Then $Q_{s+1}\in T$, which contradicts the maximality of
$s$ in (\ref{eq:maxj}) (since $s+1$ belongs to the set on the right hand side of
(\ref{eq:maxj})). This completes the proof of the Lemma.
\qed
\end{pf}
If $char\ K=0$ and $s=l_1$ then, by definition,
$\beta_s=\beta_{l_1}\ge\beta_l+\nu_l^+(f)-\nu_l(f)$. Take $q\ge l$ such that
$\delta_i=\delta$ for all $i\ge q$. Thus Lemma \ref{deltadrops} implies
that if $char\ K=0$ and $i_0\ge q$ then there exists $i>i_0$ with
$\beta_i>\beta_{i_0}+\nu_{i_0}^+(f)-\nu_{i_0}(f)$. Thus to complete the proof of the
Proposition, it remains to show (\ref{eq:infinity}) assuming that there is no $i_0$ satisfying
(\ref{eq:limbounded}).

To do that, we will define a sequence of integers $l_0,l_1,\dots$ recursively as follows. Let
$l_0=l$, where we choose $l$ sufficiently large so that $\delta_i(f)$ and $\epsilon_i(f)$
stabilize for all $i\ge l$. Let $\delta=\delta_i(f)$ and $\epsilon=\epsilon_i(f)$. By
assumption, there exists $l_1$ of the form $l+t$, $t\in\mathbb N$, such that
$\beta_{l_1}\ge\beta_l+\frac1{p^e}(\nu_l^+(f)-\nu_l(f))$. We iterate this procedure. In
other words, assume that the ordinal $l_q$ is already defined. Choose $l_{q+1}$ of the form
$l+t$, $t\in\mathbb N$, such that
\begin{equation}
\beta_{l_{q+1}}\ge\beta_{l_q}+\frac1{p^e}\left(\nu_{l_q}^+(f)-\nu_{l_q}(f)\right).\label{eq:increment}
\end{equation}
\begin{lem}\label{inc} We have
\begin{equation}
\nu_{l_1}^+(f)-\nu_{l_1}(f)\ge\nu_l^+(f)-\nu_l(f).\label{eq:step1}
\end{equation}
\end{lem}
\begin{pf} By Proposition \ref{stability},
$\nu'(a_{\delta l_1})=\nu'(a_{\delta l_0})$ and
$\nu'(a_{\epsilon l_1})=\nu'(a_{\epsilon l_0})$. Hence
$$
\begin{array}{rcl}
\nu_{l_1}^+(f)-\nu_{l_1}(f) & = & \nu'(a_{\epsilon l_1})-\nu'(a_{\delta l_1})-
(\epsilon-\delta)\beta_{l_1}\ge\\
&\ge& \nu'(a_{\epsilon l})-\nu'(a_{\delta l})-
(\epsilon-\delta)\beta_l\ge\nu_l^+(f)-\nu_l(f),
\end{array}
$$
and the Lemma is proved.
\qed
\end{pf}
We are now in the position to finish the proof of Proposition \ref{propinfty}.
Lemma \ref{inc} shows that $\nu_{l_j}^+(f)-\nu_{l_j}(f)$ is an increasing function of
$j$, so, by (\ref{eq:increment}), $\beta_{l_{j+1}}-\beta_{l_j}$ is bounded
below by an increasing function of $j$. This proves
that $\lim\limits_{q\to\infty}\beta_q=\infty$, as desired. \qed
\end{pf}
Two things remain to be accomplished in our study of infinite sequences $\{Q_{l+t}\}_{t\in\mathbb N}$ of key polynomials.
First, we must show that if $\lim\limits_{t\to\infty}\beta_{l+t}=\infty$ and $\delta=\delta_{l+t}(f)$ for $t$ sufficiently
large then $\delta$ is of the form $\delta=p^e$ for some $e\in\mathbb N$. Secondly, we must investigate the case when the
sequence $\beta_{l+t}$ is bounded and define the next key polynomial $Q_{l+\omega}$. Our main technique for dealing with
the first of these problems will be differential operators. As for the second problem, we will use Proposition \ref{diffop}
(particularly, equation (\ref{eq:close-1})). There we will not use differential operators as such, however, we will apply to
$f$ what could intuitively be termed ``differentiation of order $p^e$ with respect to $Q_i$''. We now make a digression
devoted to differential operators and their effect on key polynomials.

\section{Key polynomials and differential operators}
\label{diff}

As we saw in the previous section, the most difficult situation to handle is
one in which $t=\alpha_{i+1}=1$ in (\ref{eq:factorizl}) and (\ref{eq:alphal+1}): it is the only one which can give rise to infinite sequences of key polynomials. Then $\init_ih$ has the form
\begin{equation}
\init_ih=\init_{\nu'}d_{\delta i}(\bar Q+\init_{\nu'}z_i)^\delta=\init_iQ_{i+1}^\delta.\label{eq:initif1}
\end{equation}
This section is devoted to proving some basic results about the effect of differential operators on key polynomials, needed to study equations $h$ of the above form. Here and below, for a non-negative integer $b$, $\partial_b$ will denote the differential operator $\frac1{b!}\frac{\partial^b}{\partial x^b}$. We are interested in proving lower bounds on the quantity $\nu'(\partial_{p^b}h)$ and also in giving sufficient conditions under which $\partial_{p^b}h$ is not identically zero.

Fix an ordinal $l$ and a natural number $t$ such that
$$
\delta_{l+1}(h)=\delta_{l+2}(h)=\dots=\delta_{l+t}(h)
$$
By Proposition \ref{stability}, this implies that
\begin{equation}
\alpha_{l+2}=\dots=\alpha_{l+t}=1\label{eq:alpha=1}
\end{equation}
and that $h$ satisfies (\ref{eq:initif1}) for $l+1\le i<l+t$. Let $\delta=\delta_{l+1}(h)$. Write $\delta=p^eu$, where if
$char\ K>0$ then $p\ \not|\ u$. If $char\ K=0$, let $e=0$.

Take an ordinal $i$ having an immediate predecessor and such that the key polynomials $\mathbf Q_{i+1}$ are defined. If
$char\ K>0$, let
\begin{equation}
e_i=\min\left\{e'\ \left|\ \partial_{p^{e'}}Q_i\not\equiv0\right.\right\}.
\end{equation}
If $char\ K=0$, let $e_i=0$. Let
\begin{equation}
b=e+e_i.\label{eq:b}
\end{equation}
In the next section we will use our results on differential operators to prove that if
$\lim\limits_{t\in\mathbb N}\beta_{l+t}=\infty$ then $\delta$ is of the form $\delta=p^e$, that is, $u>1$. This will be
proved by contradiction: we will assume that $u>1$ and show that $\lim\limits_{t\in\mathbb N}\nu'(\partial_{p^b}f)=\infty$
but $\partial_{p^b}f\not\equiv0$.

Let $\mathbf Q_{i+1}^{\mbox{\boldmath$\gamma$}_{i+1}}$ be an $i$-standard monomial. One of our main tasks in this section is
to study the quantity $\nu'\left(\partial_{p^b}\mathbf Q_{i+1}^{\mbox{\boldmath$\gamma$}_{i+1}}\right)$. Since an exact
formula for $\nu'\left(\partial_{p^b}\mathbf Q_{i+1}^{\mbox{\boldmath$\gamma$}_{i+1}}\right)$ seems too complicated to
compute, we are only able to give an approximate lower bound, except under the additional assumption that
$\beta_i\gg\beta_l$ (a precise form of this inequality is (\ref{eq:betailarge}) below).

Let $b$ be any non-negative integer such that $b\ge e_i$.
\begin{prop}\label{derivative} (1) We have
\begin{equation}
\nu'\left(\mathbf Q_{i+1}^{\mbox{\boldmath$\gamma$}_{i+1}}\right)-\nu_i\left(\partial_{p^b}\mathbf
Q_{i+1}^{\mbox{\boldmath$\gamma$}_{i+1}}\right)\le
\max\left\{p^{b-e_i}\left(\beta_i-\nu_i\left(\partial_{p^{e_i}}Q_i\right)\right),p^b\beta_l\right\}.\label{eq:derivative}
\end{equation}
(2) Assume that
\begin{equation}
p^{b-e_i}\left(\beta_i-\nu_i\left(\partial_{p^{e_i}}Q_i\right)\right)>p^b\beta_l.\label{eq:betailarge}
\end{equation}
Then equality holds in (\ref{eq:derivative}) if and only if
\begin{equation}
\binom{\gamma_i}{p^{b-e_i}}\ne0.\label{eq:binom}
\end{equation}
In particular, $\partial_{p^b}\mathbf Q_{i+1}^{\mbox{\boldmath$\gamma$}_{i+1}}\not\equiv0$.

(3) Assume that both (\ref{eq:betailarge}) and (\ref{eq:binom}) hold. Then
$$
\init_i\partial_{p^b}\mathbf Q_{i+1}^{\mbox{\boldmath$\gamma$}_{i+1}}=\init_i\frac{\mathbf
Q_{i+1}^{\mbox{\boldmath$\gamma$}_{i+1}}\left(\partial_{p^{e_i}}Q_i\right)^{p^{b-e_i}}}{Q_i^{p^{b-e_i}}}.
$$
\end{prop}
\begin{pf} Direct calculation, using induction on $i$ and the fact that, by (\ref{eq:alpha=1}) and Proposition
\ref{notinvolving}, the $i$-standard monomial $\mathbf Q_{i+1}^{\mbox{\boldmath$\gamma$}_{i+1}}$ does not involve any of
$Q_{l+1}$, \dots, $Q_{i-1}$. \qed
\end{pf}
\begin{rem} The following is a well known characterization of the inequation (\ref{eq:binom}). Let
$\gamma_i=k_0+pk_1+\dots+p^qk_q$, with $k_0,\dots,k_q\in\{0,1,\dots,p-1\}$, denote the $p$-adic expansion of $\gamma_i$.
Then (\ref{eq:binom}) holds if and only if $k_{b-e_i}>0$. In particular, (\ref{eq:binom}) holds whenever $\gamma_i$ is of
the form $\gamma_i=p^{b-e_i}u$, with $p\ \not|\ u$. This is the only situation in which Proposition \ref{derivative} will
be applied in this paper.
\end{rem}
\begin{cor}\label{derivative1} Let $h$ be any element of $K[X]$, not necessarily satisfying (\ref{eq:initif1}).
(1) We have
\begin{equation}
\nu_i\left(\partial_{p^b}h\right)\ge
\nu_i(h)-\max\left\{p^{b-e_i}\left(\beta_i-\nu_i\left(\partial_{p^{b-e_i}}Q_i\right)\right),p^b\beta_l\right\}.
\label{eq:dervalue}
\end{equation}
(2) Write $\init_ih=\sum\limits_{j\in S_i}\init_i\left(d_{ji}Q_i^j\right)$. Let $S_{bi}=\left\{j\in S_i\
\left|\ \binom j{p^{b-ei}}\ne0\right.\right\}$. Assume that the inequality (\ref{eq:betailarge}) holds and that
$S_{bi}\ne\emptyset$. Then equality holds in (\ref{eq:dervalue}) and
$\init_i\partial_{p^b}h=
\sum\limits_{j\in S_{bi}}\init_i\left(d_{ji}Q_i^{j-p^{b-e_i}}\left(\partial_{p^{e_i}}Q_i\right)^{p^{b-e_i}}\right)$.
In particular, $\partial_{p^b}h\not\equiv0$.
\end{cor}
\begin{cor}\label{ppowers} Assume that $h$ satisfies (\ref{eq:initif1}). Let $b$ be as in (\ref{eq:b}). Then
$\partial_{p^b}h\not\equiv0$.
\end{cor}
\begin{pf} This is a special case of Corollary \ref{derivative1}. \qed
\end{pf}
\begin{cor} Assume that $h$ and $b$ satisfy the hypotheses of Corollary \ref{derivative1} (or, more specifically, those
of Corollary \ref{ppowers}). Then
\begin{equation}
h\notin K\left[X^{p^{b+1}}\right].\label{eq:ppowers}
\end{equation}
\end{cor}

\section{Sequences of key polynomials whose values tend to infinity}
\label{tendstoinfty}

Let the notation be as above. Let $l$ be an ordinal and assume that the above construction of key polynomials gives rise
to a sequence $\{Q_{l+t}\}_{t\in\mathbb N}$ of key polynomials such that
\begin{equation}
\lim\limits_{t\to\infty}\beta_{l+t}=\infty.\label{eq:toinfty2}
\end{equation}
Let $\delta=\delta_{l+t}(f)$ for $t$ sufficiently large. The purpose of this section is to prove
\begin{thm} The integer $\delta$ is of the form $\delta=p^e$ for some $e\in\mathbb N$.
\end{thm}
\begin{pf} We give a proof by contradiction. Suppose that (\ref{eq:toinfty2}) holds but $\delta$ is
of the form $\delta=p^eu$ with $u>1$. Let $b$ be as in (\ref{eq:b}) and let $g=\partial_{p^b}f$. The quantity $p^b\beta_l$
is independent of $t$, hence, by (\ref{eq:toinfty2}), the inequality (\ref{eq:betailarge}) holds for $t$ sufficiently
large. By Proposition \ref{stability} (2), $\init_{l+t}f$ has the form (\ref{eq:initif1}) for $i=l+t$, as $t$ runs over
$\mathbb N$. Hence $h=f$ satisfies the hypotheses of Corollary \ref{ppowers}. By Corollary \ref{ppowers}, $g\not\equiv0$.
Moreover, by Corollary \ref{derivative1} (1), we have
$\nu'(g)\ge\nu_{l+t}(g)\ge\delta\beta_{l+t}-p^e\beta_{l+t}=p^e(u-1)\beta_{l+t}$. Since $u>1$, this shows that
$\nu'(g)=\infty$, which contradicts the fact that $g$ is given by a polynomial in $x$ of degree strictly less than $n$.
\qed
\end{pf}
The following Proposition will come in useful in the remaining sections.
\begin{prop}\label{largevalue} Take an element $h$ of $L$ and an ordinal $i$ such that the key polynomials $\mathbf Q_{i+1}$
are defined. Assume that
\begin{equation}
\nu'(h)<\beta_i\label{eq:hlessbetai}
\end{equation}
and that $h$ admits an $i$-standard expansion
\begin{equation}
h=\sum\limits_{j=0}^sc_jQ_i^j,\label{eq:hstandard}
\end{equation}
such that $\nu'(c_j)\ge0$ for all $j$. Then $\nu'(h)=\nu_i(h)$.
\end{prop}
\begin{pf} By definition of standard expansion, each $c_i$ in (\ref{eq:hstandard}) is an\linebreak
$i$-standard expansion not involving $Q_i$. Then $c_j$ is a sum of monomials in $\mathbf Q_i$, which
does not vanish in $G_{\nu'}$ (\ref{eq:nonzero}), hence all the monomials appearing in $c_j$ have value at least
$\nu'(c_j)$. By (\ref{eq:hlessbetai}),
\begin{equation}
\nu'\left(c_jQ_i^j\right)=\nu_i\left(c_jQ_i^j\right)>\nu'(h)\qquad\text{ for }j>0\label{eq:greater}
\end{equation}
(\ref{eq:hstandard}) and (\ref{eq:greater}) imply that $\nu'(h)=\nu'(c_0)$. Thus $h$ is a sum of monomials in $\mathbf Q_i$
of value at least $\nu'(h)$, as desired. \qed
\end{pf}

\section{Sequences of key polynomials with bounded values in fields of positive characteristic}
\label{positive}

In this section, we assume that $char\ K=p>0$. Let $l$ be an ordinal number
and assume that the key polynomials $\mathbf Q_l\cup\{Q_{l+t}\}_{t\in\mathbb N}$ are
already defined. Moreover, assume that we are in Case 2b of \S\ref{key} (in particular,
the sequence $\{\beta_{l+t}\}_{t\in\mathbb N}$ has a upper bound $\bar\beta$ but no maximum in
$\Gamma$; this is the only case which remains to be treated to complete the definition of
the $Q_i$). By Proposition \ref{alpha1}, there exists $t_0\in\mathbb N$ such that
\begin{equation}
\alpha_{l+t}=1\ \text{ and }\ \delta_{l+t}=\delta_{l+t_0}\ \text{ for all }\ t\ge
t_0.\label{eq:deltastab}
\end{equation}
Replacing $l$ by $l+s$ for a suitable positive integer $s$, we may
assume that $\alpha_{l+t}=1$ for all strictly positive $t$.
In what follows, the index $i$ will run over the set $\{l+t\}_{t\in\mathbb N}$. As
usual, let $\delta$ denote the common value of all the $\delta_i(f)$.
\begin{prop}\label{no3a} Assume we are in Case 2b. There exist $i\in\{l+t\}_{t\in\mathbb N}$, a strictly
positive integer $e_0\le e$ and a weakly affine $i$-standard expansion
$Q_{l+\omega}$, monic of degree $p^{e_0}$ in $Q_i$, such that
\begin{equation}
\bar\beta\le\frac1{p^{e_0}}\nu\left(Q_{l+\omega}\right).\label{eq:barbeta}
\end{equation}
\end{prop}
Of course, the inequality (\ref{eq:barbeta}) is equivalent to saying that
\begin{equation}
\nu'\left(Q_{l+\omega}\right)>p^{e_0}\nu\left(Q_l+z_l+\dots+z_{l+t}\right)
\label{eq:barbeta1}
\end{equation}
for all $t\in\mathbb N$.
\begin{pf} The idea is to start with the
inequality $\nu'(f)>\nu_{l+t}(f)$ for all $t\in\mathbb N$ and to gradually construct
polynomials $g$ of the smallest possible degree satisfying
\begin{equation}
\nu'(g)>\nu_i(g)\label{eq:stricti}
\end{equation}
until we arrive at $g=Q_{l+\omega}$ satisfying the conclusion of the Proposition.

First, let $a^*$ be an $l$-standard expansion, not involving $Q_l$, such that
\begin{equation}
\init_{\nu'}(a^*a_{\delta l})=1\label{eq:adeltainverse}
\end{equation}
and let $a^*(X)$ be the representative of $a^*$ in $K[X]$ of degree less than $n$. Note that
\begin{equation}
\init_{\nu'}a_{\delta l}=\init_{\nu'}a_{\delta i}\quad\text{ for all }i\ge l\label{eq:initstable}
\end{equation}
by Proposition \ref{stability} (2).

Let $\tilde f=a^*(X)\bar f$. By Proposition \ref{stability} (2), for all $i\ge l$ we have
$$
\init_i\tilde f=\init_if=\init_{\nu'}a_{\delta i}(\bar Q_i+\init_{\nu'}z_i)^\delta,
$$
hence in view of (\ref{eq:initstable}) we have $\init_i\tilde f=(\bar Q_i+\init_{\nu'}z_i)^\delta$. In
particular,
\begin{equation}
\nu'(\tilde f)>\nu_i(\tilde f)\quad\text{ for all }i.\label{eq:nu'>nuitilde}
\end{equation}
Let $\tilde f=\sum\limits_{j=0}^{\tilde n_l}\tilde a_{ji}Q_i^j$ be
the $i$-standard expansion of $\tilde f$. We have $\init_{\nu'}\tilde a_{\delta i}=1$ for all $i$.

As noted in the previous section, since $\alpha_i=1$ for all $i$, all the $i$-standard
expansions of $\tilde f$ have the same degree $\tilde n_l$ in $Q_i$.

By Lemma \ref{inc} the quantity $\nu_i^+(\tilde f)-\nu_i(\tilde f)$ is increasing with $i$. Taking into account the fact
that $\bar\beta=\lim\limits_{i\to\infty}\beta_i$, we have, for $i$ sufficiently large,
\begin{equation}
\nu'(\tilde a_{i\delta})+\delta\bar\beta-\nu_i(\tilde f)=\delta(\bar\beta-\beta_i)<\nu_i^+(\tilde f)-\nu_i(\tilde
f).\label{eq:epsilon}
\end{equation}
By choosing $l$ sufficiently large, we may assume that (\ref{eq:epsilon}) holds for $i\ge l$.

Next, write $\tilde a_{\delta l}=1+\tilde a^\dag$ with $\nu'(\tilde a^\dag)>0$. Let $\tilde{\tilde f}=(1-\tilde
a^\dag(X))\tilde f$ and let $\tilde{\tilde f}=\sum\limits_{j=0}^{\tilde{\tilde n_l}}\tilde{\tilde a}_{jl}Q_l^j$ be the
$l$-standard expansion of $\tilde{\tilde f}$. By (\ref{eq:epsilon}) terms of the form $(1-\tilde a^\dag(X))\tilde a_{jl}$ with $j>\delta$ contribute terms of negligibly high value to $\tilde{\tilde a}_{\delta l}$. Terms $(1-\tilde
a^\dag(X))\tilde a_{jl}$ with $j<\delta$ contribute terms of value at least
$\nu'(\tilde a^\dag)+(\beta_l-\alpha_l\beta_{l-1})$ to $\tilde{\tilde a}_{\delta l}$. Thus
$\nu'(\tilde{\tilde a}^\dag)\ge\nu'(\tilde a^\dag)+\min\{\nu'(\tilde
a^\dag),\beta_l-\alpha_l\beta_{l-1}\}$, so multiplying $\tilde f$ by $(1-\tilde a^\dag(X))$ increases $\nu'(\tilde
a^\dag)$ by a fixed amount. Iterating this procedure finitely many times, we may assume that
$\nu'(\tilde a^\dag)>\delta\bar\beta>\nu_i(\tilde f)$ for all $i$. Then replacing $a_{\delta l}$ by 1 does not affect the
inequality (\ref{eq:nu'>nuitilde}), hence we may assume that $\tilde a_{\delta l}=1$.

Let
$$
\bar f=\sum\limits_{j=0}^\delta\tilde a_{jl}Q_l^j.
$$
(\ref{eq:epsilon}) implies that for all $j$, $\delta<j\le\tilde n_l$,
$$
\nu'\left(\tilde a_{jl}Q_j^l\right)\ge\nu_l^+(\tilde f)>\delta\bar\beta>\delta\beta_i=\nu_i(\tilde f).
$$
Hence $\init_i\tilde f=\init_i\bar f$; in particular, $\nu\left(\bar f\right)>\nu_i\left(\bar f\right)$ for all $i$.

The polynomial $\bar f$ is monic of degree $\delta$; the expression $\bar
f=\sum\limits_{j=0}^\delta\tilde a_{ji}Q_i^j$ is the $i$-standard expansion of $\bar f$. None of the subsequent
transformations $Q_i=Q_l+z_l+\dots+z_{i-1}$ affect the coefficient $a_{\delta l}=1$, so $a_{\delta i}=1$ for all $i$.

Write $\delta=p^eu$, as in the previous section. Let $g=\sum\limits_{j=0}^{p^e}\binom{\delta-p^e+j}j\tilde
a_{\delta-p^e+j}Q_l^j$ (roughly speaking, the reader should think of the process of constructing $g$ from $\bar f$ as
applying a differential operator of order $\delta-p^e$ with respect to $Q_{i_0}$). By construction,
\begin{equation}
\init_i\bar f=(\bar Q_i+\init_{\nu'}z_i)^\delta\qquad\text{ for }i\ge l.\label{eq:initf}
\end{equation}
On the other hand, let $w_{i-1}=z_l+\dots+z_{i-1}$. Then
\begin{equation}
\bar f=\sum\limits_{j=0}^\delta\tilde a_{jl}(Q_i-w_{i-1})^j.\label{eq:fQiw}
\end{equation}
The terms in (\ref{eq:fQiw}) with $j<\delta-p^e$ give rise to polynomials of degree strictly less than
$(\delta-p^e)\deg_xQ_l$. Thus (\ref{eq:fQiw}) can be rewritten as
\begin{eqnarray}
\bar f&=&\sum\limits_{j=0}^{p^e}\tilde a_{\delta-p^e+j,l}(Q_i-w_{i-1})^{\delta-p^e+j}+\phi\\
&=&Q_i^{\delta-p^e}\sum\limits_{j=0}^{p^e}\sum\limits_{v=j}^{p^e}(-1)^{v-j}\binom{\delta-p^e+v}{\delta-p^e+j}
\tilde a_{\delta-p^e+v,l}w_{i-1}^{v-j}Q_i^j+\psi\\
&=&Q_i^\delta+Q_i^{\delta-p^e}\sum\limits_{j=0}^{p^e-1}\sum\limits_{v=j}^{p^e}(-1)^{v-j}
\binom{\delta-p^e+v}{\delta-p^e+j}\tilde a_{\delta-p^e+v,l}w_{i-1}^{v-j}Q_i^j+\psi\label{eq:ftranslate}
\end{eqnarray}
where $\deg_x\phi,\deg_x\psi<(\delta-p^e)\deg_xQ_l$. On the other hand, we have
\begin{eqnarray}
g&=&\sum\limits_{j=0}^{p^e}\binom{\delta-p^e+j}j\tilde
a_{\delta-p^e+j,l}Q_l^j\\
&=&\sum\limits_{j=0}^{p^e}\binom{\delta-p^e+j}j\tilde a_{\delta-p^e+j,l}(Q_i-w_{i-1})^j\\
&=&\sum\limits_{j=0}^{p^e}\sum\limits_{v=j}^{p^e}(-1)^{v-j}\binom{\delta-p^e+v}v\binom vj\tilde
a_{\delta-p^e+v,l}w_{i-1}^{v-j}Q_i^j\\
&=&uQ_i^{p^e}+\sum\limits_{j=0}^{p^e-1}\sum\limits_{v=j}^{p^e}(-1)^{v-j}\binom{\delta-p^e+v}v\binom vj\tilde
a_{\delta-p^e+v,l}w_{i-1}^{v-j}Q_i^j.\label{eq:gtranslate}
\end{eqnarray}
Now, $\binom{\delta-p^e+v}v\binom vj=\binom{\delta-p^e+v}{\delta-p^e+j}\binom{\delta-p^e+j}j$ whenever $j\le v$; moreover,
\begin{eqnarray}
\binom{\delta-p^e+j}j&=&1\qquad\text{if }0\le j<p^e\\
&=&u\qquad\text{if }j=p^e.
\end{eqnarray}
Thus the double sums in (\ref{eq:ftranslate}) and (\ref{eq:gtranslate}) are identical; note also that everything in these
double sums has degree strictly less than $p^e\deg_xQ_l$. Thus rewriting the double sum as an $i$-standard expansion and
comparing (\ref{eq:gtranslate}) with (\ref{eq:initf}) shows that $\init_ig=u\bar Q_i^{p^e}+u\init_{\nu'}z_i^{p^e}=u(\bar
Q_i+\init_{\nu'}z_i)^{p^e}$; in particular, $g$ satisfies (\ref{eq:stricti}). Dividing $g$ by the non-zero integer $u$ does
not change the problem, so we may assume that $g$ is a monic polynomial in $Q_i$ of degree $p^e$. Write
\begin{equation}
g=\sum\limits_{j=0}^{p^e}c_{ji}Q_i^j.\label{eq:g}
\end{equation}
Choose $i_0\ge l$ sufficiently large so that
\begin{equation}
\beta_{i_0}-\alpha_l\beta_{l-1}>p^e(\bar\beta-\beta_{i_0}).\label{eq:i0close}
\end{equation}
\begin{rem}\label{replace} Assume that there exist $i\ge i_0$ and $j$, $1\le j<p^e$, such that
$\nu'(c_{ji})+j\bar\beta>2p^e\bar\beta-p^e\beta_i$. Then for any $i'>i$ we have $\init_{i'}(g-c_{ji}Q_i^j)=\init_{i'}g$; in
particular, $\nu_{i'}(g-c_{ji}Q_i^j)<\nu'(g-c_{ji}Q_i^j)$. Thus we are free to replace $g$ by $g-c_{ji}Q_i^j$.
\end{rem}
Assume that there exist $j\in\{1,...,p^e-1\}$ and $i_1\ge i_0$ such that $c_{ji_1}\ne0$ and
\begin{equation}
p^e\bar\beta<\nu(c_{ji_1})+j\bar\beta<(p^e+1)\beta_{i_1}-\alpha_l\beta_{l-1}.\label{eq:gap}
\end{equation}
Take the greatest such $j$.
\begin{lem}\label{dichotomy1} (1) We have $\nu'(c_{ji})+j\bar\beta>p^e\bar\beta$ for all $i\ge
i_1$.

(2) The element $\init_{\nu'}c_{ji}$ is constant for all $i\ge i_1$.

(3) There exists $i_2\ge i_1$ such that for all $i\ge i_2$ we have
$$
\nu_i\left(g-c_{ji_2}Q_{i_2}^j\right)<\nu'\left(g-c_{ji_2}Q_{i_2}^j\right).
$$
\end{lem}
\begin{pf} (2) follows the maximality of $j$ and the inequalities (\ref{eq:i0close}) and (\ref{eq:gap}):
$\init_{\nu'}c_{ji_2}$ cannot affected by any subsequent coordinate changes of the form
$Q_i=Q_{i_1}+z_{i_1}+\dots+z_{i-1}$. (1) follows immediately from (2).

By (1) and (2), taking $i_2$ sufficiently large, we can ensure that
$$
\nu'\left(c_{ji_2}Q_{i_2}^j\right)>p^e\bar\beta.
$$
Since $p^e\bar\beta>p^e\beta_i=\nu_i(g)$, we have
$$
\nu_i\left(g-c_{ji_2}Q_{i_2}^j\right)=\nu_i(g)<
\min\left\{\nu'(g),\nu'\left(c_{ji_2}Q_{i_2}^j\right)\right\}\le\nu'\left(g-c_{ji_2}Q_{i_2}^j\right)
$$
for all $i\in i_2+\mathbb N$, and (3) is proved. This completes the proof of Lemma \ref{dichotomy1}. \qed
\end{pf}
If there exists $j\in\{1,\dots,p^e-1\}$ satisfying the hypotheses of Lemma \ref{dichotomy1}, replace $g$ by
$g-c_{ji_2}Q_{i_2}^j$; Lemma \ref{dichotomy1} (3) says that strict inequality (\ref{eq:stricti}) is satisfied with $g$
replaced by $g-c_{i_2j}Q_{i_2}^j$. This procedure strictly decreases the integer $j$ appearing in Lemma \ref{dichotomy1}.
Hence after finitely many repetitions of this procedure we obtain a polynomial $g$ such that there do not exist $j$ and
$i_1$ satisfying (\ref{eq:gap}). By the second inequality in (\ref{eq:gap}), the non-existence of such $j$ and $i_1$ is
preserved as we pass from $i$ to $i+1$; hence, after finitely many steps we may assume that no $j$ and $i_1$ satisfying
(\ref{eq:gap}) exist. We will make this assumption from now on.
\begin{rem} Now, by the same reasoning as in Lemma \ref{dichotomy1}, the sets
$$
\mathcal{S}:=\left\{\ j\in\{1,...,p^e\}\left|\ c_{ji}\ne0\text{ and }\nu(c_{ji})=(p^e-j)\bar\beta\right.\right\}
$$
and $\left\{\left.\init_{\nu'}c_{ji} \right|\ j\in\mathcal{S}\right\}$ are independent of $i$ for $i\ge i_0$.
\end{rem}
\begin{lem}\label{dichotomy} Consider an index $j\in\{1,...,p^e-1\}$ and an ordinal $i\ge i_0$ of the form $i=i_0+t$,
$t\in\mathbb N$, as above. Assume that $c_{ji}\ne0$. We have
\begin{equation}
\nu'(c_{ji})+j\bar\beta\ge p^e\bar\beta\label{eq:aboveline}
\end{equation}
and $j$ is a power of $p$ whenever equality holds in (\ref{eq:aboveline}).
\end{lem}
\begin{pf} We give a proof by contradiction. Assume that for a certain $i_1\ge i_0$ there exists $j\in\{1,...,p^e-1\}$ such
that $c_{ji_1}\ne0$, and either
\begin{equation}
\nu(c_{ji_1})<(p^e-j)\bar\beta\label{eq:below}
\end{equation}
or $j$ is not a $p$-power (or both). Let $j(g)$ denote the greatest such $j$. Let $j=j(g)$.
Then the element $\init_{\nu'}c_{ji_1}$ is not affected by the subsequent coordinate
changes $Q_{i_1}=Q_i-z_{i_1}-\dots-z_{i-1}$, so $\init_{\nu'}c_{ji}=\init_{\nu'}c_{ji_1}$ for all $i\ge i_1$.

First assume that (\ref{eq:below}) holds. (\ref{eq:below}) can be rewritten as $\nu(c_{ji})+j\bar\beta<p^e\bar\beta$. Now,
taking $i$ sufficiently large, the difference $\bar\beta-\beta_i$ can be made arbitrarily small, so
$\nu(c_{ji})+j\beta_i<p^e\beta_i$. This inequality shows that
$$
\nu\left(c_{ji}Q_i^j\right)<\nu\left(Q_i^{p^e}\right),
$$
so $\init_ig$ does not contain the monomial $Q_i^{p^e}$, which is a contradiction.

From now on assume that
\begin{equation}
\nu'(c_{ji})+j\bar\beta=p^e\bar\beta\qquad\text{ for all }i\ge i_1.\label{eq:online}
\end{equation}
Then, by definition of $j$, $j$ is not a $p$-power. Write $j=p^{e'}u'$ and
$Q_{i+1}=Q_i+z_i$. Then the $(i+1)$-standard expansion of $g$ contains a
monomial of value $\nu'\left(c_{ji}z_i^{j-p^{e'}}Q_{i+1}^{p^{e'}}\right)$. We have
$$
\nu'\left(c_{ji}z_i^{j-p^{e'}}\right)+p^{e'}\bar\beta=
\nu'(c_{ji})+\left(j-p^{e'}\right)\beta_i+p^{e'}\bar\beta<\nu'(c_{ji})+j\bar\beta=
p^e\bar\beta.
$$
Thus the appearance of a monomial of value $\nu'\left(c_{ji}z_i^{j-p^{e'}}Q_{i+1}^{p^{e'}}\right)$ in the
standard expansion of $g$ contradicts (\ref{eq:online}) with $i$ replaced by $i+1$. This completes the proof of Lemma
\ref{dichotomy}. \qed
\end{pf}
If $Q_{l+\omega}=g$ satisfies the conclusion of Proposition \ref{no3a} there is nothing more
to prove. Otherwise, by Lemma \ref{dichotomy} and since no $j$ and $i_1$ satisfy
(\ref{eq:gap}), there exist $j\in\{1,\dots,p^e-1\}$ and $i_1\ge i_0$ such that for all $i\ge i_1$ we have
\begin{equation}
\nu(c_{ji})+j\bar\beta>(p^e+1)\beta_i-\alpha_l\beta_{l-1}\ge(p^e+1)\beta_{i_1}-\alpha_l\beta_{l-1}>p^e\bar\beta.
\label{eq:gap1}
\end{equation}
Let $\mathcal{A}$ denote the set of all such $j$. Replace $g$ by $g-\sum\limits_{j\in\mathcal{A}}c_{ji_1}Q_{i_1}^j$. Remark
\ref{replace} says that strict inequality (\ref{eq:stricti}) is satisfied for this new $g$.
In this way, we obtain a polynomial $g$ such that $Q_{l+\omega}=g$ satisfies the conclusion of Proposition
\ref{no3a}. This completes the proof of Proposition \ref{no3a}. \qed
\end{pf}
\begin{rem} We are not claiming that the property that $g$ is a weakly affine expansion in $Q_{i_1}$ is preserved when
we pass from $i_1$ to some other ordinal $i>i_1$. However, the above results show that for any $i\ge i_1$ of the form
$i=l+t$, $t\in\mathbb N$, $g$ is a sum of a weakly affine expansion in $Q_i$ all of whose monomials $c_{ji}Q_i^j$ lie
on the critical line $\nu'(c_{ji})=(p^e-j)\bar\beta$ and another standard expansion of degree strictly less than $p^e$ in
$Q_i$, all of whose monomials have value greater than or equal to $(p^e+1)\beta_{i_1}-\alpha_l\beta_{l-1}>p^e\bar\beta$.
\end{rem}
We define $Q_{l+\omega}$ to be a weakly affine standard expansion satisfying the conclusion of Proposition \ref{no3a},
which minimizes the integer $e_0$ (so that $\alpha_{l+\omega}=p^{e_0}$). This completes the definition of the $Q_i$.

Let $\theta_{l+\omega}(f)=\frac\delta{\alpha_{l+\omega}}$. It is easy to see, by the same argument as in Lemma \ref{nui},
that the Newton polygon $\Delta_{l+\omega}(f)$ contains a vertex $(\nu'(a_{\theta_{l+\omega}(f)}),\theta_{l+\omega}(f))$,
and that this vertex lies above the pivotal vertex $(\nu'(a_{\delta_{l+\omega}(f)}),\delta_{l+\omega}(f))$. The vertex
$(\nu'(a_{\theta_{l+\omega}(f)}),\theta_{l+\omega}(f))$ will be called the \textbf{characteristic vertex} of
$\Delta_{l+\omega}(f)$. The notion of characteristic vertex will be used in \S\ref{descr} when we study the totality of
extensions of $\nu$ to $L$. It is important that the characteristic vertex is determined by $\mathbf Q_{l+\omega+1}$ and
$\mbox{\boldmath$\beta$}_{l+\omega}$; it does not depend on $\beta_{l+\omega}$.
\begin{rem} By construction, we have $\alpha_{l+\omega}=p^{e_0}\ge p$. Then the fact
that $\deg_xQ_{l+\omega}\le n$ and Proposition \ref{bound} show that the situation considered
in this section can arise at most $[\log_pn]$ times, so the set $\mathbf Q:=\{Q_i\}$ thus
defined has order type of at most $[\log_pn]\omega+t$, where $t\in\mathbb N$.
\end{rem}

\section{Proof that $\{Q_i\}$ is a complete set of key polynomials}
\label{keyproperty}

This section is devoted to proving
\begin{thm}\label{mainthm} The well ordered set $\mathbf Q:=\{Q_i\}$ defined in the previous
sections is a complete set of key polynomials. In other words, for any element $\beta\in\Gamma_+$ the
corresponding $\nu'$-ideal $\mathbf P'_\beta$ is generated by all the monomials in the $Q_i$
of value $\beta$ or higher. In particular, we have $G_{\nu'}=G_\nu[\init_{\nu'}\mathbf Q]^*$.
\end{thm}
\begin{cor} The valuation $\nu'$ is completely determined by the data $\mathbf Q,\{\beta_i\}$.
\end{cor}
\begin{pf} Let $\lambda$ be the ordinal number which represents the order type of the set
$\mathbf Q$, so that $\mathbf Q=\mathbf Q_\lambda$. Let $l$ denote the smallest ordinal such that
$0\le l<\lambda$ and $\alpha_i=1$ whenever $l<i<\lambda$ (note, in particular, that if
$\lambda$ admits an immediate predecessor and $\alpha_{\lambda-1}>1$ then $l=\lambda-1$; at
the other end of the spectrum is the possibility that $\alpha_i=1$ for all $i<\lambda$ and
$l=0$). To prove the Theorem, it is sufficient to show that for every positive $\beta\in\Gamma$ and every $h\in L$
such that $\nu'(h)=\beta$, $h$ belongs to the ideal generated by all the monomials $c\mathbf Q^{\mbox{\boldmath$\gamma$}}$
such that $\nu'\left(c\mathbf Q^{\mbox{\boldmath$\gamma$}}\right)\ge\beta$.

Take any element $h\in L$. Without loss of generality, we may assume that, writing $h=\sum\limits_{j=0}^sd_jx^j$, we have
$\nu(d_j)\ge0$ for all $j$ (otherwise, multiply $h$ by a suitable element of $K$).
\begin{claim} There exists $i<\lambda$ of the form $i=l+t$, $t\in\mathbb N$, such that
\begin{equation}
\beta_i>\nu'(h).\label{eq:toinfty}
\end{equation}
\end{claim}
\begin{pf} There are two possibilities: either $\lambda$ has an immediate
predecessor or it does not. By construction, for any $i$ such that $l<i<\lambda$ we have $i=l+t$ for some
$t\in\mathbb N$. The ordinal $\lambda$ admits an immediate predecessor if and only if $\lambda=l+t$ for
some $t\in\mathbb N$ and does not admit an immediate predecessor if and only if $\lambda>l+t$
for all $t\in\mathbb N$. If $\lambda$ has an immediate predecessor then $Q_{\lambda-1}=f(x)=0$,
so $\nu'(Q_{\lambda-1})=\infty>\nu'(h)$. If $\lambda$ does not have an immediate predecessor then by construction
$\lim\limits_{t\to\infty}\beta_{l+t}=\infty$, so there exists $i=l+t$, $t\in\mathbb N$ such that (\ref{eq:toinfty}) holds.
The Claim is proved.\qed
\end{pf}
Now, Lemma \ref{largevalue} says that $\nu_i(h)=\nu'(h)$. This means, by definition, that $h$ can be written as a sum of
monomials in $\mathbf Q_{i+1}$ of value at least $\nu'(h)$, hence it belongs to the ideal generated by all such monomials.
This completes the proof.
\qed
\end{pf}

\section{A description of the algorithm.}
\label{descr}

Let $K\hookrightarrow L$ be a finite separable field extension and
$\nu:K^*\rightarrow\Gamma$ a valuation of $K$. In this section we
describe an algorithm for constructing all the possible extensions
$\nu'$ of $\nu$ to $L$.  Pick and fix a generator $x$ of $L$ over
$K$ once and for all. Let $f=\sum\limits_{i=0}^na_ix^i$ denote the
minimal polynomial of $x$ over $K$.

First, we reduce the problem to the case $rk\ \nu=1$.
Let $r=rk\ \nu$. Write $\nu$ as a composition of $r$ rank 1 valuations:
$\nu=\nu_1\circ\dots\circ\nu_r$, where $\nu_1$ is the valuation of $K$, centered at the
smallest non-zero prime ideal of $R_\nu$. Assume the problem is already solved for rank 1
valuations. Then any extension $\nu'$ of $\nu$ to $L$ is of the form
$\nu'=\nu'_1\circ\dots\circ\nu'_r$, where $\nu'_1$ is an extension of $\nu_1$ to $L$, $\nu'_2$
is an extension of $\nu_2$ to $k_{\nu'_1}$, and so on. The valuation $\nu'_i$ is an extension
of the valuation $\nu_i$ of the field $k_{\nu_{i-1}}$ to its algebraic extension
$k_{\nu'_{i-1}}$. Thus, it is sufficient to solve the problem in the case $rk\ \nu=1$.

From now on, we assume that $rk\ \nu=1$.

\noindent\textbf{Step 1.1 of the algorithm.} Choose an element $\beta_1\in\Gamma_+$ which determines a side
of $\Delta(f)$ and put $\nu'(x)=\beta_1$.

\noindent\textbf{Step 1.2 of the algorithm.} Let
\begin{equation}
\init_1f=v\prod\limits_{j=1}^sg_j^{\gamma_j}\label{eq:factorization}
\end{equation}
be the factorization of $\init_1f$ into (monic) irreducible factors in
$G_\nu[\bar x]$. Since for any extension $\nu'$ of $\nu$ we have $\init_1f(\init_{\nu'}x)=0$, one of the irreducible
factors in (\ref{eq:factorization}) is the minimal polynomial of $\init_{\nu'}x$ over $G_\nu$. Choose one of the
irreducible factors in (\ref{eq:factorization}) (other than $\bar x$), say $g_1$. Write
$$
g_1=\sum\limits_{i=0}^{\alpha_1}\bar b_i\bar x^i,
$$
where $\bar b_{\alpha_1}=1$. For each $i$, $0\le i\le\alpha_1$, let $b_i$ be a
representative of $\bar b_i$ in $R_\nu$ (that is, an element of
$R_\nu$ such that $\init_\nu b_i=\bar b_i$). Put $Q_1=x$ and
$Q_2=\sum\limits_{i=1}^{\alpha_1}b_ix^i$.

Assume, inductively, that key polynomials $Q_1$, \dots, $Q_l$ and positive
integers $\alpha_1$, \dots, $\alpha_{l-1}$ are already constructed for a certain ordinal $l$, where $l<\omega$ if $char\
K=0$ and $l<([\log_pn]+1)\omega$ if $char\ K=p>0$.

Assume, inductively, that for each $i$, $1\le i<l$, the $(i+1)$-st key
polynomial $Q_{i+1}$ admits an $i$-th standard expansion of the form
\begin{equation}
Q_{i+1}=Q_i^{\alpha_i}+\sum\limits_{j=0}^{i-1}\left(\sum\limits_\gamma
c_{ji\gamma}\mathbf Q_j^\gamma\right)Q_i^j,\label{eq:standardalg}
\end{equation}
where each of $c_{ji\gamma}\mathbf Q_j^\gamma$ is
an $i$-standard monomial. Assume that the standard expansions (\ref{eq:standardalg}) satisfy all the conditions described
in \S\ref{key}.

Write $f=\sum\limits_{j=0}^{n_l}a_{jl}Q_l^j$, where each $a_{jl}$ is a homogeneous
$l$-standard expansion not involving $Q_l$. The next two steps
of the algorithm are a generalization of the first two steps, with 1 replaced by $l$.

\noindent\textbf{Step l.1 of the algorithm.} If $l$ does not have an immediate predecessor (that is, $l$ is of the form
$l=l_0+\omega$), let $\bar\beta_l=\sup\{\beta_{l_0+t}\}_{t\in\mathbb N}$.
Choose an element $\beta_l$ which determines a side $A_l$ of $\Delta_l(f)$ and satisfies the following condition:

\textbf{Condition (*).} If $l$ has an immediate predecessor then $\beta_l>\alpha_l\beta_{l-1}$; if $l$
does not have an immediate predecessor then $\beta_l>\alpha_l\bar\beta_l$.

Put $\nu'(Q_l)=\beta_l$.
\begin{rem} We know from \S\ref{key} and \S\ref{positive} that Condition (*) must hold for any extension
$\nu'$ of $\nu$; it is a consequence of the proof of Proposition \ref{stability} (see (\ref{eq:charpivot})) that the
pivotal vertex of $\Delta_l(f)$ lies below its characteristic vertex. Conversely, if $\mathbf Q_{l+1}$ and
$\mbox{\boldmath$\beta$}_l$ are given, any side of $\Delta_l(f)$ lying below its characteristic vertex can be chosen to be
the characteristic side; the choice of $\beta_l$ which determines such a characteristic side will automatically
satisfy Condition (*).
\end{rem}
\noindent\textbf{Step l.2 of the algorithm.} By Proposition \ref{value} the value
$\nu'(a_{jl})$, where $1\le i\le n_l$, is completely determined by the $l$-standard
expansion of $a_{il}$; in particular, it is completely determined at this stage of the
algorithm. Similarly, $\init_lf:=\sum\limits_{(\nu'(a_{il}),i)\in A_l}\init_{\nu'}a_{il}\bar Q_l^i$
is a well-defined element of $G_\nu\left[\init_{\nu'}Q_1,\ldots,\init_{\nu'}Q_{l-1}\right]\left[\bar
Q_l\right]$. Let $\init_lf=v_l\prod\limits_{j=1}^{t_l}g_{jl}^{\gamma_{jl}}$ be the
factorization of $\init_lf$ into (monic) irreducible factors in
$G_\nu\left[\init_{\nu'}Q_1,\ldots,\init_{\nu'}Q_{l-1}\right]\left[\bar Q_l\right]$. Choose one
of these factors (other than $\bar Q_l$), say $g_{1l}$ (then $g_{1l}$ will be the minimal polynomial of $\init_{\nu'}Q_l$
over $G_\nu\left[\init_{\nu'}Q_1,\ldots,\init_{\nu'}Q_{l-1}\right]$ for the valuation $\nu'$
we are about to construct). Let $\alpha_{l+1}=\deg_{\bar Q_l}g_{1l}$. Write
\begin{equation}
g_{1l}=\bar Q_l^{\alpha_{l+1}}+\sum\limits_{j=0}^{\alpha_{l+1}-1}
\left(\sum\limits_{\mbox{\boldmath$\gamma$}_l}\bar
c_{l+1,j\mbox{\boldmath$\gamma$}_l}\init_{\nu'}\mathbf Q_l^{\mbox{\boldmath$\gamma$}_l}\right)
\bar Q_l^j,\label{eq:ingraded1}
\end{equation}
If $t_l>1$ or $\alpha_{l+1}>1$, define \textbf{the $(l+1)$-st key polynomial} of $\nu'$ to be a lifting
$$
Q_{l+1}=Q_l^{\alpha_{l+1}}+\sum\limits_{j=0}^{\alpha_{l+1}-1}
\left(\sum\limits_{\mbox{\boldmath$\gamma$}_l}
c_{l+1,j\mbox{\boldmath$\gamma$}_l}\mathbf Q_l^{\mbox{\boldmath$\gamma$}_l}\right)
Q_l^j
$$
(\ref{eq:ingraded1}) to $L$. If $t_l=\alpha_{l+1}=1$, the $(l+1)$-st key
polynomial $Q_{l+1}$ will also be a lifting of (\ref{eq:ingraded1}) to $L$, but we require
it to satisfy additional conditions, as in \S\ref{key}. Let $\delta_l(f)$ be defined as in
(\ref{eq:delta}). Define the next key polynomials $Q_{l+1}$, $Q_{l+2}$, \dots, as in \S\ref{key}. More precisely, we
define finitely many polynomials $Q_{l+1}$, \dots, $Q_s$ if either $char\ K=0$ or $char\ K=p>0$ and we are in Case 1 of
\S\ref{key}.

In Case 2, there exists an infinite sequence $z_l,z_{l+1},\dots$ of homogeneous standard
expansions in $\mathbf Q_l$, not involving $Q_l$, such that the sequence
$\{\nu'(Q_l+z_l+\dots+z_{l+t})\}_{t\in\mathbb N}$ is strictly increasing; pick and fix one such sequence. Define
$Q_{l+t}=Q_l+z_l+z_{l+1}+\dots+z_{l+t-1}$ for $t\in\mathbb N$. For each key polynomial $Q_i$, write
$f=\sum\limits_{j=0}^{n_l}a_{ji}Q_i^j$ and consider the corresponding Newton polygon
$\Delta_i(f)$. By definition of $\delta_i(f)$, the Newton polygon $\Delta_i(f)$ contains a
vertex $(\nu'(a_{\delta_i(f)i}),\delta_i(f))$. Since $\delta_i(f)=\delta_{i+1}(f)$, the characteristic side $A_i$ of
$\Delta_i(f)$ is uniquely determined, that is, there exists a unique element
$\beta_i\in\Gamma\cup\{\infty\}$ such that $\beta_i\ge\beta_{i'}$ for all $i'<i$, $\beta_i$ determines a
side $A_i$ of $\Delta_i(f)$ and $(\nu'(a_{\delta i}),\delta_i(f))$ is the leftmost endpoint of $A_i$.
This defines an infinite sequence $\{Q_{l+t}\}_{t\in\mathbb N}$ of key polynomials, such that for each $i=l+t$,
$t\in\mathbb N$, we have $\init_if=\init_{\nu'}a_{\delta i}(\bar Q_i+\init_{\nu'}z_i)^{\delta_l(f)}$.

Let $\bar\beta=\lim\limits_{t\to\infty}\nu'(Q_l+z_l+\dots+z_{l+t})$. By definition, we are in Case 2a if $\bar\beta=\infty$
and in Case 2b if $\bar\beta<\infty$.

In Case 2b, define the next key polynomial to be a polynomial $Q_{l+\omega}$ satisfying the conclusion of Proposition
\ref{no3a}. Note that in all the cases both the slope of the characteristic side $L_i$ and the irreducible factor
of $\init_if$ which is the minimal polynomial of $\init_{\nu'}Q_i$ over $G_\nu[\init_{\nu'}\mathbf Q_i]^*$
are uniquely determined.

The algorithm stops if one of the following occurs: either $Q_i=0$
or
$$
\sup\limits_i\{\beta_i\}=\infty,
$$
where $\beta_i$ ranges over the values of key polynomials defined so far. In both cases,
the valuation $\nu'$ is completely determined by the data $\{Q_i,\beta_i\}$.

This completes our construction of the extensions $\nu'$. Note that \textbf{every} choice described in the algorithm above
leads to an extension $\nu'$. Indeed, such a choice defines, in particular, the well ordered set $\{\nu_i\}_{i\in\Lambda}$
of valuations of $K[X]$ and their graded algebras; whenever $i<i'$, we have a natural homomorphism of graded algebras
$G_{\nu_i}\rightarrow G_{\nu_{i'}}$. The proof of Theorem \ref{mainthm} applies verbatim to show that for each $h\in L$,
the value $\nu_i(h(X))$ stabilizes for $i$ sufficiently large. Setting $\nu'(h)$ to be that stable value of $\nu_i(h(X))$
defines a valuation $\nu'$ of $L$.
\begin{cor}\label{unique} The extension $\nu'$ is unique if and only if, for each $i$ in the
above algorithm, the following two conditions hold:

(1) The $i$-th Newton polygon $\Delta_i(f)$ has only one face $L_i$ (other than the two axes).

(2) The corresponding initial form $\init_if$ does not have two distinct irreducible factors
(in other words, $\init_if$ is a power of an irreducible polynomial).
\end{cor}
The next Corollary is valid for valuations of arbitrary rank (and not only for those of rank 1).
\begin{cor} Assume that $\init_{\nu'}x$ has degree $n$ over $G_\nu$. Then $\nu$ admits a unique
extension $\nu'$ to $L$.
\end{cor}
\begin{pf} By writing $\nu$ as a composition of several rank 1
valuations, it is sufficient to prove the Corollary under the assumption
$rk\ \nu=1$. Now, the hypotheses imply that (1) and (2) of Corollary
\ref{unique} hold for $i=1$. Moreover, we may take $f=Q_2$, so the algorithm
consists of only one step, and the Corollary follows. \qed
\end{pf}
We end this paper with a discussion of the well known formula
\begin{equation}
\sum\limits_{j=1}^tf_je_jd_j=n,\label{eq:ramification}
\end{equation}
where
$\{\nu'_1,\dots,\nu'_t\}$ is the set of all the extensions of $\nu$ to $L$, $f_j$ is the index of the value group of $\nu$
inside the value group of $\nu'_j$, $e_j$ is the degree of the residue field extension $k_\nu\rightarrow k_{\nu'_j}$ and
$d_j$ is the defect of $\nu'_j$. Of course, for each $j$ we have $f_je_j=[G_{\nu'_j}:G_\nu]$. We associate to the above
algorithm the following finite, oriented, weighted tree $U$. The set of vertices of $U$ is partially ordered. In each
vertex, we have a key polynomial $Q_l$ appearing at some step in one of the branches of the above algorithm. The important
data associated to this vertex is the data $\delta_l(f)$, as well as the data $\mathbf Q_l$ of all the key polynomials
preceding $Q_l$ in the given branch of the algorithm. The set of vertices has a unique minimal element and the key polynomial
associated to this minimal vertex is $x=Q_1$. Each vertex is adjacent to exactly one vertex smaller than itself and,
possibly, to finitely many vertices greater than itself. Let us denote each vertex by the key polynomial $Q_l$ associated
to it. Not every key polynomial will be associated to a vertex of $U$. If $l$ admits an immediate predecessor then the
unique vertex adjacent to $Q_l$, preceding $Q_l$, is $Q_{l-1}$. Consider a vertex $Q_l$. We will now describe all the
vertices following $Q_l$. There are two possibilities:

(a) There is a unique $\beta_l$ satisfying Condition (*) and Case 2b of \S\ref{key} holds in the definition of $Q_l$.

(b) Condition (a) does not hold.

In case (a), the unique vertex following $Q_l$ is $Q_{l+\omega}$. In case (b), $Q_l$ is followed by all the possible
key polynomials $Q_{l+1}$, appearing in the above algorithm.

This information determines the tree $U$ completely. It is obvious that $U$ is finite.
\begin{prop}\label{branching} Fix a vertex $Q_l$ of $U$. Assume that Case (b) holds for $Q_l$ and let $Q_{l+1}^{(1)}$,
\dots, $Q_{l+1}^{(s)}$ denote all the vertices of $U$, adjacent to $Q_l$ and following it. Let $\theta_l(f)$,
$\theta_{l+1}^{(1)}(f)$, \dots, $\theta_{l+1}^{(s)}(f)$, and $\alpha_{l+1}^{(1)}$, \dots, $\alpha_{l+1}^{(s)}$ denote the
numerical characters corresponding to the $s$ resulting branches of the above algorithm. Then
\begin{equation}
\sum\limits_{j=1}^s\alpha_{l+1}^{(j)}\theta_{l+1}^{(j)}(f)=\theta_l(f).\label{eq:branching}
\end{equation}
\end{prop}
\begin{pf} Let $A_l^{(1)}$, \dots, $A_l^{(t)}$ denote the sides of $\Delta_l(f)$ lying below the characteristic
vertex $(\nu'(a_{\theta_l(f)l}),\theta_l(f))$. For $1\le j\le t$, let $\beta_l^{(j)}$ denote the element of $\Gamma_+$
which determines the side $A_l^{(j)}$ and let $\init_l^{(j)}f$ denote the corresponding initial form of $f$. By
construction, renumbering the vertices $Q_{l+1}^{(j)}$, if necessary, we can find indices $1\le s_1<\dots<s_t=s$ such that
the factorization of $\init_l^{(j)}f$ into irreducible factors has the form $\init_l^{(j)}f=\bar
Q^{u_j}\prod\limits_{q=s_{j-1}+1}^{s_j}\init_l^{(j)}(Q_{l+1}^{(q)})^{\theta_{l+1}^{(q)}}$, where the exponent $u_j$ may or
may not be zero (cf. Lemma \ref{factors} and (\ref{eq:initQi+1})). Since $\theta_l(f)$ equals the sum of the heights of
sides of $\Delta_l(f)$ lying below the characteristic vertex, we have
$\theta_l=\sum\limits_{j=1}^t\deg_{\bar
Q}\prod\limits_{q=s_{j-1}+1}^{s_j}\init_l^{(j)}\left(Q_{l+1}^{(q)}\right)^{\theta_{l+1}^{(q)}}$. Recalling that $\deg_{\bar
Q}\init_lQ_{l+1}^{(q)}=\alpha_{l+1}^{(q)}$ completes the proof of the Proposition. \qed
\end{pf}
If Case (a) holds for $Q_l$ then the pivotal vertex of $\Delta_l$ is uniquely determined and coincides with the
characteristic vertex. There is only one choice for the key polynomial $Q_{l+\omega}$ and, by definition,
\begin{equation}
\theta_l=\alpha_{l+\omega}\theta_{l+\omega}.\label{eq:branching1}
\end{equation}
Thus, the analogue of the formula (\ref{eq:branching}) holds also in the Case (a).

For each vertex $Q_l$ of $U$, let $\alpha_l(Q_l)$ denote the integer $\alpha_l$ corresponding to $Q_l$ in the above
algorithm, and similarly for $\theta_l(Q_l)$. Fix a vertex $Q_l$ of $U$ and consider a subtree $U'\subset U$, having the
following properties:

(1) $Q_l$ is the unique minimal element of $U'$.

(2) For each vertex $Q_i$ of $U'$, if $U'$ contains one vertex immediately following $Q_i$ then it contains all of them.

Let $\{Q_{l_1},\dots,Q_{l_t}\}$ be the set of maximal elements among the vertices of $U'$.
\begin{cor} We have
\begin{equation}
\theta_l=\sum\limits_{j=1}^t\left(\prod_{Q_{l'}\le Q_{l_j}}\alpha_{l'}(Q_{l'})\right)
\theta_{l_j}(Q_{l_j}).\label{eq:theta}
\end{equation}
\end{cor}
\begin{pf} This follows immediately from Proposition \ref{branching} and equation (\ref{eq:branching1}) by
induction on the size of $U'$. \qed
\end{pf}
Let $\{\nu'_1,\dots,\nu'_s\}$ be the set of all the extensions of $\nu$ to $L$ and take $U'=U$ in the above Corollary.
Let $\{Q_{l_1},\dots,Q_{l_s}\}$ be the set of maximal elements among the vertices of $U$. For each $j\in\{1,\dots,s\}$,
consider the following partition of the set of all vertices $Q_l$ of $U$ such that $Q_l\le Q_{l_j}$. We will say that such
a $Q_l$ belongs to the set $D_j$ if Case (a) holds for the vertex immediately preceding $Q_l$, and belongs to the set $E_j$
otherwise. Noting that $\theta_1=n$, we can now rewrite (\ref{eq:theta}) as
$n=\sum\limits_{j=1}^t\left(\prod\limits_{Q_{l'}\in D_j}\alpha_{l'}(Q_{l'})\right)\left(\prod\limits_{Q_{l'}\in
E_j}\alpha_{l'}(Q_{l'})\right)\theta_{l_j}(Q_{l_j})$.

Now, if $Q_{l'}\in E_j$ then the graded algebra extension
$$G_\nu[\init_{\nu'}Q_{l'}]\hookrightarrow G_\nu[\init_{\nu'}Q_{l'}][\init_{\nu'}Q_{l'}]$$ has degree $\alpha_{l'}(Q_{l'})$.
We can now interpret the formula (\ref{eq:ramification}) by observing that $\left(\prod\limits_{Q_{l'}\in
E_j}\alpha_{l'}(Q_{l'})\right)$ equals the degree of the graded algebra extension
$$
[G_{\nu'_j}:G_\nu]=e_jf_j,
$$
whereas the quantity $\left(\prod\limits_{Q_{l'}\in D_j}\alpha_{l'}(Q_{l'})\right)\theta_{l_j}(Q_{l_j})$ is nothing but the
defect of the extension $\nu'_j$.

We refer the reader to Michel Vaqui\'e's paper \cite{V4} for a detailed treatment of defect.







\end{document}